\mag=1200
\input amstex
\documentstyle{amsppt}
\document
\NoBlackBoxes
\TagsOnRight
\tolerance=2000

\topmatter

\title
Towards the Intersection Theory on Hurwitz Spaces
\endtitle
\author
M.~E.~Kazaryan, S.~K.~Lando
\endauthor
\date
01.04.2004
\enddate
\address
Steklov Mathematical Institute RAS and
\hfil\break\indent
Independent University of Moscow,
\hfil\break\indent
Institute for System Research RAS and
\hfil\break\indent
Independent University of Moscow,
\endaddress

\abstract
Moduli spaces of algebraic curves and closely related to
them Hurwitz spaces, that is, spaces of meromorphic functions
on the curves, arise naturally in numerous problems of algebraic
geometry and mathematical physics, especially in relationship with 
the string theory and Gromov--Witten invariants. In particular,
the classical Hurwitz problem about enumeration of topologically 
distinct ramified coverings of the sphere with prescribed ramification
type reduces to the study of geometry and topology of these spaces. 
The cohomology rings of such spaces are complicated
even in the simplest cases of rational curves and functions.
However, the cohomology classes that are the most important from the 
point of view of applications (namely, the classes Poincar\'e dual
to the strata of functions with given singularities)
can be expressed in terms of relatively simple ``basic''
classes (which are, in a sense, tautological). 
The aim of the present paper is to identify these basic classes,
to describe relations among them, and to find expressions for
the strata in terms of these classes. Our approach is based on
R.~Thom's theory of universal polynomials of singularities,
which has been extended to the case of multisingularities
by the first author. Although the general Hurwitz problem 
still remains open, our approach allows one to achieve
a significant progress in its solution, as well as in the 
understanding of the geometry and topology of Hurwitz spaces.

\endabstract
\endtopmatter
\rightheadtext{On the Intersection Theory on Hurwitz Spaces}

\define\ch{\operatorname{ch}}
\define\Aut{\operatorname{Aut}}
\define\Hom{\operatorname{Hom}}
\def\dts{\dots,}
\define\td{\operatorname{td}}

\footnote""{The research was supported by the RFBR grants
(no.~04-01-00762 for the first author and no.~02-01-22004
for the second one).
}

\head
\S\,1. Introduction
\endhead

\subhead 1.1. Hurwitz's problem\endsubhead
In~\cite{8} A.~Hurwitz suggested the following problem:
enumerate isomorphism classes of ramified coverings
of the $2$-sphere by a surface of genus~$g$ 
having prescribed ramification points in the target sphere
and prescribed ramification types over each ramification
point. (The formal statement of the problem is given below.)
During the last century the problem attracted researchers'
attention several times, however this direction of mathematical
research became a central one only recently. The explosion of
interest to the problem is due first of all
to the applications it has found in modern mathematical physics,
namely, in the string theory and the theory of Gromov--Witten invariants.

Up to now, the Hurwitz problem does not have a satisfactory solution:
known enumerative formulas lead to transparent answers 
only in some partial cases. In~\cite{1},~\cite{2} 
a relationship between ramified coverings of the sphere 
and intersection theory on moduli spaces of functions
on algebraic curves has been established. The very existence of
such connection revives hopes for obtaining simple general
answers. In \cite{1},~\cite{2} the Hurwitz problem was reduced
to a problem concerning the geometry of moduli spaces;
the latter is of its own interest. Because of the construction
in~\cite{2} it is closely related to the already classical
problem about the geometry of moduli spaces of algebraic curves.

In the present paper we analyze the part of the geometry
of Hurwitz spaces directly related to the Hurwitz problem.

\subhead
1.2. Hurwitz spaces and cohomology classes of strata
\endsubhead
We shall consider ramified coverings of the $2$-sphere
by surfaces of genus~$g$; let~$n$ denote the degree
of the coverings. In addition to topological coverings, we
shall also consider meromorphic functions of degree~$n$
on algebraic curves of genus~$g$. Of course, from
the topological point of view such a function is a ramified
covering. Denote by~$\Cal H_{g,n}$ the space of such functions
possessing the following properties:

--~all the poles of the functions are simple, i.e., 
the meromorphic function has~$n$ poles of order~$1$;
moreover, we assume that the poles of each function are numbered;

--~the sum of the critical values of the function is~$0$.

According to~\cite{2} this space is a smooth complex orbifold
(or even a smooth variety provided that
$g=0$ or~$n$ is sufficiently large). It is fibered over
the moduli space~$\Cal M_{g,n}$ of complex curves of genus~$g$
with~$n$ marked points: to each function one can associate its
underlying curve equipped with the~$n$ marked poles.

The space~$\Cal H_{g,n}$ possesses a completion denoted 
by~$\overline{\Cal H}_{g,n}$ and consisting of~{\it stable\/}
meromorphic functions \cite{2},~\cite{14}. Its {\it boundary\/}
$\overline{\Cal H}_{g,n}\setminus \Cal H_{g,n}$ 
consists of stable meromorphic functions with a singular 
underlying curve, the only admissible singularities being nodes
(points of double selfintersection).
The projection $\Cal H_{g,n}\to\Cal M_{g,n}$ extends to a projection
$\overline{\Cal H}_{g,n}\to\overline{\Cal M}_{g,n}$ 
of the completed space to the moduli space of
{\it stable\/} curves with marked points. The fibers
of this projection are vector spaces since a linear combination
of meromorphic functions with poles of order at most one
at the marked points is a function of the same kind. 
Note, however, that generally speaking this projection
is not a vector bundle since the dimension of the fiber
can vary from one point of the base to another. 
The fiberwise projectivization $P\overline{\Cal H}_{g,n}$ 
of the completed Hurwitz space is a compact complex orbifold.
This will be our main moduli space, and all the other required
spaces will be constructed using it. By some abuse of language,
we speak below about varieties and subvarieties, having
in mind that our spaces are, in fact, orbifolds and suborbifolds
(i.e., locally they look like quotient spaces of the 
complex disc modulo a finite group action).

Since the variety~$P\overline{\Cal H}_{g,n}$
is a projectivization, it carries a natural second cohomology class,
namely, the first Chern class of the tautological sheaf,
which we denote by
$\psi=\psi_{g,n}=c_1(\Cal O(1))\in H^2(P\overline{\Cal H}_{g,n})$.

On the other hand, this variety contains subvarieties
consisting of functions with degenerate ramification.
By the Riemann--Hurwitz formula, a generic meromorphic function
of degree~$n$ on a curve of genus~$g$ has $2n+2g-2$
points of nondegenerate ramification (both in the source
and in the target). The functions with fewer ramification points
in the target sphere form the
{\it discriminant\/} in the space of functions.
One can assign to each ramification point in the target sphere
a partition of~$n$ which is the (unordered) set of multiplicities
of the preimages of this point. It will be more convenient
for us, however, to use the {\it reduced partitions},
that is, the sets of multiplicities of the preimages decreased by one.
We denote the closure in~$P\overline{\Cal H}_{g,n}$ 
of the set functions having ramification of prescribed type 
by~$\sigma_{\alpha_1,\dots,\alpha_c}$; the subscript
here consists of reduced partitions over the points
of degenerate ramification. These subvarieties are called
the {\it strata\/} of the discriminant. Thus, 
$\sigma_{2^1}\subset P\overline{\Cal H}_{g,n}$ 
denotes the {\it caustic}, i.e., the stratum consisting
of functions with two ramification points in the preimage glued together,
and~$\sigma_{1^2}$ is the {\it Maxwell stratum}
whose open part consists of functions taking the same value
at two distinct critical points. The caustic and the Maxwell stratum 
are the only strata of (complex) codimension one.
Similarly, the stratum $\sigma_{2^1,1^2}$, of codimension~$2$,
is the closure of the set of functions with two degenerate
critical values, with a ramification point of multiplicity~$3$
over one of them, and two simple ramification points over the other one, 
and so on. The number of nondegenerate critical points
(corresponding to the reduced partition~$1^1$)
for a generic function in a given stratum can be easily calculated
by means of the Riemann--Hurwitz formula. For the sake of simplicity
of notation we do not include these nondegenerate reduced partitions,
as well as the values~$g$ and~$n$, in the notation for the strata.

Each stratum is a complex subvariety of pure
codimension in~$P\overline{\Cal H}_{g,n}$, 
whence it determines, by Poincar\'e duality, a homogeneous
class in the cohomology ring
$$
H^*(P\overline{\Cal H}_{g,n})=H^*(P\overline{\Cal H}_{g,n},\Bbb Q).
$$
(Everywhere below we are interested in cohomology with
rational coefficients, and we omit this basic ring
in the notation for cohomology.) The degree of this element
coincides with the (real) codimension of the stratum.
The product of this element with
the class~$\psi$ taken to the complementary power
is simply a rational number
which, and this is the result of the geometrical approach
suggested in~\cite{2}, is closely related to the
number of ramified coverings we are interested in.
Before describing this connection let us give the rigorous
statement of the Hurwitz problem.

Two ramified coverings  $f_1\:C_1\to S^2$, \ $f_2\:C_2\to S^2$
of the $2$-sphere are said to be {\it isomorphic\/} 
if there is a homeomorphism $h\:C_1\to C_2$ such that $f_1=f_2\circ h$.
Clearly, isomorphic ramified coverings have coinciding
ramification points in the target sphere and coinciding ramification
types over these points, as well as the degrees~$n_1,n_2$
and genera~$g_1,g_2$ of the covering surfaces,
$n_1=n_2=n$, \ $g_1=g_2=g$. A ramification point in the target
sphere is said to be {\it nondegenerate}, or {\it simple},
if it has~$n-1$ geometrically distinct preimages one of which
is a ramification point of multiplicity~$2$,
and the other $n-2$ points are smooth points of the covering.
As it was mentioned above, we assign to each
ramification point in the target the reduced partition~$\alpha$
consisting of the multiplicities of its singular preimages
decreased by one. The set of all reduced partitions~$\alpha$ 
for a given ramified covering is called the {\it ramification data}.
The {\it Hurwitz number\/} 
$h_{\alpha_1,\alpha_2,\dots}$ is defined by the equation
$$
h_{\alpha_1,\alpha_2,\dots}=\sum_f\frac1{\bigl|\Aut(f)\bigr|}\,,
$$
where the summation is carried over all isomorphism classes
of ramified coverings~$f$ of the sphere of degree~$n$
by a surface of genus~$g$, with reduced partitions~$\alpha_1,\alpha_2,\dots$
over the degenerate ramification points, and
$|\Aut(f)|$ denotes the number of elements in the automorphism group
of such ramified covering (this group is always finite).
Of course, the Hurwitz numbers depend also on~$g$ and~$n$,
but we do not include these parameters in the notation.

The {\it Hurwitz problem\/} consists in counting Hurwitz numbers,
and it is related to the geometry of Hurwitz spaces
by the following theorem, which has been explicitly
formulated, for the case $g=0$, in~\cite{16}
and which follows immediately from the results in~\cite{2}.

\proclaim{Theorem 1.1}
We have
$$
h_{\alpha_1,\alpha_2,\dots}=
\frac{\bigl|\Aut(\alpha_1,\alpha_2,\dots)\bigr|}{n!}\,
\langle\sigma_{\alpha_1,\alpha_2,\dots},\psi^d\rangle,
$$
where $\bigl|\Aut(\alpha_1,\alpha_2,\dots)\bigr|$
denotes the order of the automorphism group of the
set of reduced partitions, that is, the product
of factorials of the numbers of pairwise coinciding partitions,
and the angle brackets denote the multiplication
in the cohomology of $P\overline{\Cal H}_{g,n}$.
\endproclaim

Theorem~1.1 reduces the Hurwitz problem to a problem
about the cohomology of the Hurwitz spaces. 
In spite of the visible complication caused by replacing
a group-combinatorial problem with an algebro-geometric one,
this geometric interpretation not only clarifies the known
computational results (see \cite{8}, \cite{6}, \cite{15},~\cite{2}),
but also leads to new ones (see, e.g.,~\cite{16}).

Generally speaking, the cohomology rings of moduli spaces
of stable curves and closely related to them Hurwitz spaces are
very complicated, even in the case of rational curves and functions
($g=0$). Our experience shows, however, that all natural
cohomology classes related to singularities (as well as multisingularities
and multimultisingularities --- these are exactly the classes
playing the central role in the Hurwitz problem) can be expressed in terms
of a small relatively simple set of ``basic'' (tautological, in a sense)
classes. Our goal is to identify these basic classes, describe 
relations between them, and find expressions for strata in terms 
of these classes.

\goodbreak
\subhead 1.3. Description of main results\endsubhead
We start with describing results concerning Hurwitz numbers.

\proclaim{Theorem 1.2}
For rational
{\rm(}i.e., related to the case $g=0${\rm)} Hurwitz numbers the
following formulas are true{\rm:}
$$
\align
h_{2^1,2^1}&=\frac34 (27 n^2-137 n+180)\frac{(2 n-6)!}{(n-3)!}\,n^{n-6},
\\
h_{2^1,1^2}&=3(3 n^2-15 n+20)\frac{(2 n-6)!}{(n-4)!}\, n^{n-6},
\\
h_{1^2,1^2}&=2 (2 n^3-16 n^2+43 n-40)\frac{(2 n-6)!}{(n-4)!}\,n^{n-6}.
\endalign
$$
\endproclaim

These formulas describe some multimultisingularities, that is,
the case where the ramification is nondegenerate over 
{\it several\/} points in the target sphere. 
For the case of multisingularities, where the ramification
is degenerate over a {\it single\/} point, the rational Hurwitz numbers
are given by the well-known Hurwitz formula. The first formula in
Theorem~1.2 was obtained by D.~Zvonkine in~\cite{21},
the other two, as far as we know, are new.

To prove these results, we deduce some relations in the cohomology
of Hurwitz spaces on the classes represented by the discriminant strata
and the class~$\psi$. We shall also make an extensive use of the
class $\delta\in H^2(P\overline{\Cal H}_{g,n})$
represented by functions on singular curves. The class~$\delta$ 
is very close to the pullback of the boundary class
in~$\overline{\Cal M}_{g,n}$ with respect to the projection
$\overline{\Cal H}_{g,n}\to\overline{\Cal M}_{g,n}$ 
(the difference of the two is the divisor of functions
having a single pole on a rational irreducible component 
attached to other components at a single point). 
The value of the latter in the geometry of moduli spaces
of stable curves is well-known. For example, the index
of maximal selfintersection of this class is nothing but the
{\it Weil--Petersson\/} volume of the moduli space
(see~\cite{17}). The essential difference between Hurwitz
spaces and moduli spaces of curves is arising
due to existence of another natural class of $2$-cohomology, namely,
the class~$\psi$.

Note that all the results presented in the paper are valid
in the Chow ring of
$P\overline{\Cal H}_{g,n}$, not only in the cohomology. The only reason
why we speak about cohomology is that the cohomological computations
are sufficient for the expected computational applications.

As far as we know, the first relations for strata in
the cohomology rings of Hurwitz spaces
were obtained in~\cite{16}. They express the caustic and the 
Maxwell stratum in the spaces of rational functions
($g=0$) in terms of the boundary stratum and the class~$\psi$:
$$
\aligned
\sigma_{2^1}&=6(n-1)\psi-3\delta,
\\
\sigma_{1^2}&=2(n-1)(n-6)\psi+4\delta.
\endaligned
\tag 1.1
$$

We present below a new derivation of these formulas.
The classes whose complex codimension is greater than one also
are subject to numerous relations which are described
in more detail in the main body of the paper,
since their description requires additional definitions.

\subhead 1.4. Degree\endsubhead
Theorem~1.1 reduces the process of calculating Hurwitz numbers
to the computation of the intersection index of a stratum
with the complementary power of the class~$\psi$.
Such class~$\psi$ arises naturally on any orbit space
of an action of the group~$\Bbb C^*$; 
it behaves functorially under equivariant mappings of varieties
endowed with actions of this group. That is why
the intersection index of a subvariety
with the complementary degree of this class 
deserves the name of {\it degree\/} of the subvariety.

Let $A$~be a complex variety and let the multiplicative group~$\Bbb C^*$
of nonzero complex numbers act on~$A$ without fixed points.
Let~$B$ be the orbit space of this action and suppose~$B$ is compact.
Denote by~$\Psi=c_1\bigl(O(1)\bigr)\in H^2(B)$ the tautological
class of this action. For an arbitrary class
$\beta\in H^*(B)$, the {\it degree\/}~$\deg\beta$ is
the result of the pairing
$$
\deg\beta=\int_B\frac{\beta}{1-\Psi}=
\int_B\beta(1+\Psi+\Psi^2+\Psi^3+\dots).
$$
In other words, the degree of a homogeneous class $\beta$
is the degree of its intersection with the class~$\Psi$
taken to the power equal to the codimension of~$\beta$.

\example{Example 1.3}
The weighted projective space $\Bbb CP_{w_0,\dots,w_n}$
is the quotient space of the complement to the origin
in the vector space $\Bbb C^{n+1}$ endowed with the following action
of the group~$\Bbb C^*$:
$$
\lambda\:(x_0,\dots,x_n)\mapsto(\lambda^{w_0}x_0,\dots,\lambda^{w_n}x_n),
\qquad\lambda\in\Bbb C^*.
$$
The integral exponents $w_i$ are the {\it weights\/} 
of the action. Each holomorphic linear representation
of~$\Bbb C^*$ has such a form. The degree of the weighted
projective space $\Bbb CP_{w_0,\dots,w_n}$
with the weights $w_0,\dots,w_n$ is $\bigl(\prod w_i\bigr)^{-1}$.
\endexample

The definitions immediately imply the functoriality of the degree.

\proclaim{Proposition 1.4}
Suppose the group~$\Bbb C^*$ acts without fixed points
on complex varieties~$A_1,A_2$ of the same dimension, and
suppose the orbit spaces $B_1,B_2$ are compact. Then
the degree of an equivariant mapping $f\:A_1\to A_2$ 
is independent of~$f$ and can be determined from the equality
$$
\deg B_1=\deg f\cdot\deg B_2.
$$
\endproclaim

Theorem~1.1 is an immediate corollary of this Proposition
applied to the mapping
$L_{g,n}\:\overline{\Cal H}_{g,n}\to\Bbb C^{2g+2n-3}$
taking a meromorphic function to the unordered set
of its critical values (recall that according to our
agreement the sum of the critical values is zero).

\subhead
1.5. Singularities, relative Chern classes, and universal
polynomials
\endsubhead
Equations~\thetag{1.1} were obtained in~\cite{16}
by computing the intersection indices of the classes~$\sigma_{2^1}$,
$\sigma_{1^2}$ with a basic set of classes of complementary dimension
(i.e., of complex dimension one, that is, curves).
However as the codimension of the strata grows this way of arguing becomes
cumbersome: both the choice of a basic set of classes
of complementary dimension and the computation of
the intersection index with such a class are complicated problems,
without known approaches to their effective solutions.
We present a new approach to the search for relations
between cohomology classes. This approach is based on Thom's theory
of universal polynomials for singularities~\cite{20}
extended to the case of multisingularities in~\cite{10}--\cite{13}.
Let us explain briefly the basics of this theory.

Let $F\:M\to N$ be a generic holomorphic mapping of
compact complex varieties (we assume, for definiteness,
that the varieties have coinciding dimensions: although the theory
works for arbitrary dimensions, we are going to apply it only in
the case of coinciding dimensions). The {\it total Chern class\/}
of~$F$ is the ratio
$$
\split
c(F)&=\frac{c(F^*TN)}{c(TM)}=1+(c_1(F^*(TN))-c_1(TM))+(c_2(F^*TN)
\\
&\qquad
-c_1(F^*TN)c_1(TM)-c_2(TM)+c_1^2(TM))+\ldots\in H^*(M);
\endsplit
$$
here the function $c(\,\cdot\,)$ whose argument is a vector
bundle denotes, as usual, the total Chern class of this bundle.
The homogeneous components of the total Chern class
are called the {\it Chern classes\/} of~$F$
and denoted by $c_i(F)\in H^{2i}(M)$:
$$
\aligned
c_1(F)&=c_1(F^*TN)-c_1(TM),
\\
c_2(F)&=c_2(F^*TN)-c_1(F^*TN)c_1(TM)-c_2(TM)+c_1^2(TM),
\endaligned
$$
and so on. A mapping~$F$ has a singularity at a given point
if its differential~$dF$ is degenerate at this point.
Singularities are split into various types.
The closure of the set of points where~$F$ has a singularity of
a given type determines a cohomology class in~$H^*(M)$.
Thom's theory implies that this class can be expressed
as a polynomial in the Chern classes of~$F$, and the coefficients
of the polynomial depend only on the type of the singularity,
not on the varieties~$M$ and~$N$ or the mapping~$f$.

The notion of both a generic mapping and a type of a singularity
requires a more precise definition which can be found, e.g., 
in~\cite{13}. As the simplest example let us consider
the singularity~$A_1$: a mapping has singularity~$A_1$
at a given point if the rank of the differential of~$F$
at this point is one less than the maximal possible,
and the restriction of~$F$ to the kernel of the differential
is nondegenerate. For a general mapping the set of points
where it has a singularity of type~$A_1$ is a subvariety
of complex codimension one. The cohomology class
of this closure coincides with~$c_1(F)=c_1(F^*TN)-c_1(TM)$,
which can be easily shown by considering the pullback
to~$M$ of a nonzero differential form of the highest degree on~$N$.
Already in the simplest cases this expression for the class~$A_1$
leads to nontrivial results. Suppose, for example, that 
the dimensions of both varieties is one, i.e.,
we deal with ramified coverings of degree, say, $n$
of a complex curve by a complex curve. Then we obtain
the Riemann--Hurwitz formula for the number of ramification points
of a generic ramified covering:
$$
c_1(F)=n(2-2g_N)-(2-2g_M)=2(n-1)-2(ng_N-g_M),
$$
where $g_M$ and~$g_N$ are, respectively, the genus of the covering
and of the covered surface.

\goodbreak
A similar statement for {\it multisingularities\/} asserts
that the closure in~$N$ of the set of points such that~$F$
has at their preimages singularities of prescribed types
can be expressed as a universal polynomial in the pushforwards
$F_*((c_1(F))^{k_1}(c_2(F))^{k_2}\dots)\in H^*(N)$
of monomials in Chern classes of~$F$. We also extend this
approach to the case of {\it multimultisingularities},
that is, subvarieties in moduli spaces of functions
consisting of functions having prescribed sets of multisingularities.

We apply the methods of the theory of universal polynomials
to the universal mapping over the Hurwitz space
which we are going to define. Denote 
by~$\overline{\Cal U}_{g,n}$ the {\it universal curve\/}
over~$\overline{\Cal H}_{g,n}$, i.e., a complex orbifold together
with a mapping $\overline{\Cal U}_{g,n}\to\overline{\Cal H}_{g,n}$
whose fiber over a point $f\in\overline{\Cal H}_{g,n}$ 
is the domain of definition of~$f$ (or, more precisely, the 
quotient of this domain modulo the action of the automorphism
group of the function~$f$). One more variety is just the 
direct product 
$\Bbb CP^1\times\overline{\Cal H}_{g,n}$, 
with a coordinate chosen in the first factor. The {\it universal
mapping\/} over the Hurwitz space takes a point of the universal 
curve~$\overline{\Cal
U}_{g,n}$ to the value of the corresponding function at this point.
It is fibered over the Hurwitz space.

\remark{Remark 1.5}
As it is well known, the universal curve over the moduli space
of stable curves~$\overline{\Cal M}_{g,n}$ 
coincides with the moduli space
$\overline{\Cal M}_{g,n+1}$ of stable curves with~$n+1$ 
marked points. Indeed, a point of a fiber of the universal curve
can be treated as the $(n+1)$~st marked point along the fiber. 
A similar construction works for Hurwitz spaces as well.
\endremark

\medskip

The multiplicative group~$\Bbb C^*$ of nonzero complex numbers
acts on each of the three spaces~$\overline{\Cal H}_{g,n}$, \
$\overline{\Cal U}_{g,n}$ and~$\Bbb CP^1\times\overline{\Cal H}_{g,n}$
by multiplying a function by a constant. On~$\overline{\Cal U}_{g,n}$
this action is fiberwise trivial, while on the fiber~$\Bbb CP^1$ 
of the direct product the group~$\Bbb C^*$ 
acts as multiplication by constants (this action is well defined
since we have chosen a coordinate in the fiber). The universal mapping
is equivariant with respect to the above actions, and 
the strata in the Hurwitz space are invariant with respect to the action.
After restricting both the universal curve and the direct product
to the complement to the set of fixed points of the action on~$\overline{\Cal
H}_{g,n}$ and taking the quotients of the resulting spaces 
modulo the group actions, we arrive at the triangle of spaces and mappings
which will be the main object of our study. Note that since the
action of~$\Bbb C^*$ on the fibers of the direct product 
$\Bbb CP^1\times\overline{\Cal H}_{g,n}$ is nontrivial,
its quotient does not coincide with the direct product~$\Bbb CP^1\times P
\overline{\Cal H}_{g,n}$; instead, it is isomorphic to the projectivization
of the vector bundle 
$\Bbb C\oplus\Cal O(1)$ of rank two over $P\overline{\Cal H}_{g,n}$.

Introduce the simplified notation: $B_{g,n}$ for the projectivized Hurwitz
space $P\overline{\Cal H}_{g,n}$, \ $X_{g,n}$
for the quotient of the universal curve (which, in its own turn, is
the universal curve over the projectivized Hurwitz space),
$Y_{g,n}$ for the quotient of the direct product,
$p_{g,n}$ (respectively,
$q_{g,n}$) for the projection of $X_{g,n}$ (respectively,
$Y_{g,n}$) to~$B_{g,n}$ and~$f_{g,n}$ for the quotient
universal mapping $f_{g,n}\:X_{g,n}\to Y_{g,n}$. In other words,
we consider the commutative triangle
$$
\alignedat2
&\quad \;X_{g,n}&\overset{f_{g,n}}\to\longrightarrow\quad&\quad Y_{g,n}
\\
&\qquad{\ssize{p_{g,n}}}\!\!\!\searrow &\swarrow&\ssize{q_{g,n}}
\\
&&B_{g,n}&
\endalignedat
\tag 1.2
$$

We start with showing that the relative Chern classes
of the universal mapping can be expressed in terms of few
``basic'' classes in~$H^*(X_{g,n})$, and then describe
the strata in the space~$B_{g,n}$ in terms of the 
pushforwards of these basic classes. This proves to be realizable
because the mapping~$f_{g,n}$ has not so many types of singularities,
and the types are independent of~$g$ and~$n$, which leads to
universal formulas.

Note that since the Hurwitz space~$B_{g,n}$ is the moduli space
of stable meromorphic functions of degree~$n$ 
on curves of genus~$g$ and~$X_{g,n}$ is the universal curve on~$B_{g,n}$, 
any family of stable meromorphic functions of degree~$n$
on curves of genus~$g$ is induced from the projection
$p_{g,n}\:X_{g,n}\to B_{g,n}$. In particular, linear relations 
for cohomology classes valid in~$X_{g,n}$
and~$B_{g,n}$ remain valid in arbitrary family.

\subhead 1.6. The structure of the paper\endsubhead
In~\S\,2 we analyze local singularities of the
universal mapping and introduce the basic cohomology classes.
In~\S\,3 we show that the relative Chern classes
of the universal mapping indeed can be expressed in
terms of the basic classes. Then we deduce formulas
for some strata in~$B_{g,n}$ in terms of the pushforwards
of polynomials in basic classes. In~\S\,4 
the degrees of the strata are computed.

\medskip

The authors express their gratitude to B.~Lass, D.~Zagier and 
D.~Zvonkine for useful discussions. The paper was completed
during the second author's stay at the Max-Planck Institut
f\"ur Mathematik, Bonn, Germany, in January--February 2004. 

\head
\S\,2. Singularities of fibered mappings with one-dimensional fibers
\endhead

The relative Chern classes of the universal mapping~$f_{g,n}$
are closely related to its singularities. Since the
mapping~$f_{g,n}$ is fibered over the base~$B_{g,n}$
and the fibers are one-dimensional, these singularities
are relatively simple. The classification of simple
isolated singularities with one-dimensional fibers~\cite{5}
includes singularities of types~$A_k$
and~$I_{k,l}$. The mappings~$f_{g,n}$ acquire such singularities.
However in addition to isolated singularities these mappings
also have (except for the simplest cases
$g=0$ and $n=2$ or $n=3$) nonisolated ones. The latter
arise if a fiber which is the domain of the function
contains an irreducible component where the function is constant.
Note that such a function is not necessarily unstable: 
if the genus of the irreducible component is at least one, or
it is rational but contains at least three nodes, then 
the function can be constant on such a component, and yet stable. 
Of course, our approach to the study of the geometry 
of Hurwitz spaces requires a complete investigation of
nonisolated singularities, which we do not present here,
whence we obtain only partial results. We plan to discuss
nonisolated singularities in a separate paper. 

\subhead 2.1. Local singularities\endsubhead
Consider the commutative triangular diagram of spaces
and holomorphic mappings of the form
$$
\alignedat2
&\Bbb C^{m+1} &\overset{f}\to\longrightarrow\quad&\Bbb C\times\Bbb C^m
\\
&\qquad\;{\ssize{p}}\!\searrow&\qquad&\swarrow\ssize{q}\quad
\\
&&\Bbb C^m&
\endalignedat
\tag 2.1
$$
where $q$ is the projection of the direct product to the second
factor. We assume that the fibers of~$p$ are one-dimensional
(although not necessarily smooth). Restricting~$f$ to the fibers
of~$p$ we obtain a family of functions on (noncompact)
curves. The natural ``right'' equivalence of such diagrams
admits arbitrary coordinate changes in the base~$\Bbb C^m$,
as well as fiberwise coordinate changes in~$\Bbb C^{m+1}$
fibered over the chosen coordinate transformation in the base.
We extend this equivalence group allowing for 
multiplicating a function by a constant depending holomorphically
on the point of the base. This extension does not affect 
the classification of singularities, but leads to a
richer theory of characteristic classes.

Since the mapping~$q$ is a projection, and~$p$ is the composition 
of~$f$ and~$q$, the singularities of the diagram are
totally determined by those of~$f$. A point~$x\in\Bbb C^{m+1}$
is a {\it local singularity\/} of~$f$ if the rank of the
differential of~$f$ at~$x$ is less than~$m+1$. Let us make
the following two assumptions:

-- the mapping~$f$ has only generic singularities, i.e.,
the complement to the set of such singularities
in the space of jets of diagrams is a closed nowhere dense subset;

--  the projection $p\:\Bbb C^{m+1}\to\Bbb C^m$ 
has no singularities that are more complicated than the Morse
folding $(x,y)\mapsto xy$.

For the universal mapping we are interested in,  the second assumption
is satisfied, while the first one is not: the restriction of
the function to a fiber can prove to be constant on some 
irreducible component of the fiber. As a result, the formulas
we obtain are valid ``up to some classes supported on the subvariety
of nonisolated singularities''. In the case of rational functions
the codimension of the subvariety of nonisolated singularities
is~$3$, which allows us to obtain exact formulas for the strata
up to codimension two.

Under the above assumptions the mapping~$f$ can have singularities
of the following types:

(1) the type~$A_k$, \ $k=1,2,\dots$, is realized at those smooth points
of the fibers of~$p$, where the restriction of~$f$ to the fiber
can be reduced to the form~$x\mapsto x^{k+1}$
in an appropriate local coordinate~$x$ along the fiber;

(2) the type~$I_{k,l}$, \ $k\geq l\geq 1$, is realized at those
singular points of the fibers of~$p$, where the restriction
of~$f$ to the plane with the local coordinates~$x,y$,
such that~$p$ has the form $(x,y)\mapsto xy$ looks like
$f\:(x,y)\mapsto (x^k+y^l,xy)$.

At those points one can choose coordinates on the base~$\Bbb C^m$
such that~$f$ can be represented as the direct sum of the identity mapping
and miniversal unfoldings which, for the cases~$A_k$
and~$I_{k,l}$, look as follows:

(1) \  $A_k\:(x,b_2,\dots,b_{k})
\mapsto(x^{k+1}+b_2x^{k-1}+\dots+b_kx,b_2,\dots,b_k)$;

(2) \ $I_{k,l}\:(x,y,a_1,\dots,a_{k-1},b_1,\dots,b_{l-1})\mapsto
(x^k+a_1x^{k-1}+\dots+a_{k-1}x+y^l+b_1y^{l-1}+
\dots+b_{l-1}y,xy,a_1,\dots,a_{k-1},b_1,\dots,b_{l-1})$.

We interpret these unfoldings as globally given
polynomial mappings of complex vector spaces
of coinciding dimensions. In the case~$A_k$
(respectively, $I_{k,l}$) such an unfolding is the universal
mapping of the universal curve over the moduli
space of rational functions with a single pole of order~$k+1$
(respectively, with two poles, of orders~$k$ and~$l$)~\cite{14}.
More precisely, this is not the moduli space, but some
finite covering of the moduli space, of degree~$k+1$ 
in the case~$A_k$ and of degree~$kl$ in the case~$I_{k,l}$.
Therefore, in the above coordinates the multiplication of
a function by a constant produces the multiplication of
a rational function of the type~$A_k$ (respectively, $I_{k,l}$)
by~$\lambda^{k+1}\in\Bbb C^*$ (respectively, $\lambda^{kl}$)
which, after the substitution $x\mapsto\lambda x$,
defines the action of the group~$\Bbb C^*$ on the unfoldings
according to the following formulas:

(1) for the singularity~$A_k$
$$
\lambda\:(x,b_2,\dots,b_k)\mapsto(\lambda x,\lambda^2b_2,\dots,
\lambda^kb_k),\qquad \lambda\in\Bbb C^*;
$$

(2) for the singularity~$I_{k,l}$
$$
\split
\lambda&\:(x,y,a_1,\dots,a_{k-1},b_1,\dots,b_{l-1})
\\
&\qquad\mapsto
\bigl(\lambda^lx,\lambda^ky,\lambda^la_1,\dots,\lambda^{l(k-1)}a_{k-1},
\lambda^kb_1,\dots,\lambda^{k(l-1)}b_{l-1}\bigr),\qquad\lambda\in\Bbb C^*.
\endsplit
$$

\subhead 2.2. Characteristic classes of singularities\endsubhead
Suppose that diagram~\thetag{2.1} has a singularity of the type~$A_k$
at the origin. Without loss of generality one can suppose that
$k=m+1$, and~$f$ is exactly the above miniversal unfolding 
of the singularity~$A_k$. The space~$\Bbb C^k$ in the upper left
corner of diagram~\thetag{2.1} contains subvarieties 
corresponding to all singularities~$A_i$ for
$i\leq  k$. These subvarieties are invariant with respect to
the action of the group~$\Bbb C^*$, therefore,
to the singularity~$A_i$ a subvariety~$[A_i]$ 
in the weighted projective space
$\Bbb CP^{(k)}=\Bbb CP_{1,2,\dots,k}$ obtained by factorization
of the punctured space~$\Bbb C^k\setminus\{0\}$ modulo the action
of the group~$\Bbb C^*$ described above is assigned, 
whence a cohomology class~$[A_i]\in H^{2i}(\Bbb CP^{(k)})$. 
Denote by~$\Sigma\in H^2(\Bbb CP^{(k)})$
the characteristic class of the singularity~$A_1$, 
i.e., the first Chern class of~$f$, 
and by $\Psi\in H^2(\Bbb CP^{(k)})$ the characteristic class 
of the action of~$\Bbb C^*$.

\proclaim{Theorem 2.1}
The characteristic classes of the singularities~$[A_i]$
in the unfoldings of the singularities~$A_k$ 
can be expressed in terms of $\Sigma$ and~$\Psi$
according to the formula
$$
[A_i]=P_i(\Sigma,\Psi)=
\Sigma(2\Sigma-\Psi)(3\Sigma-2\Psi)\dots\bigl(i\Sigma-(i-1)\Psi\bigr).
$$
In particular, the expressions for these classes
in terms of $\Psi$ and $\Sigma$ are independent of the order~$k$
of the perturbed singularity~$A_k$.
\endproclaim

\demo{Proof}
Denote by~$t\in H^2(\Bbb CP^{(k)})$ the cohomology class such that
$\Psi=(k+1)t$. Then~$t$ is a generator of the cohomology ring
$H^*(\Bbb CP^{(k)})$. Because of quasihomogenity,
we have $\Sigma=kt$. The expression for the 
class~$[A_i]$ in terms of~$t$ and~$k$ must have the form
$$
[A_i]=c_i\cdot(k)_it^i
$$
for some constant~$c_i$. Here the Pohgammer symbol $(k)_i$
denotes the following polynomial in~$k$ of degree~$i$:
$$
(k)_i=k(k-1)(k-2)\dots(k-i+1).
$$

Indeed, Thom's theorem together with Corollary~3.4 below
imply that the class~$[A_i]$ admits some universal
representation as a homogeneous polynomial in~$\Psi$ 
and~$\Sigma$. Substituting into this polynomial 
the expressions for~$\Psi$ and~$\Sigma$ in terms of~$t$ and~$k$, 
we conclude that the class~$[A_i]$
is proportional to~$t^i$ and the proportionality coefficient
is a polynomial in~$k$ of degree~$i$ vanishing for all~$k<i$. 
For $i=k$ we have~$[A_i]=i!\,t^i$ which yields $c_i=1$.
Now the expression for~$A[i]$ in terms of~$\Psi$ and~$\Sigma$
results from the substitution  $t=\Psi-\Sigma$, \
$k=\Sigma/(\Psi-\Sigma)$. The theorem is proved.
\enddemo

The miniversal unfolding of the singularity~$I_{k,l}$ 
contains both the strata corresponding to the singularities~$A_i$ 
for~$i<k+l$ and the strata corresponding to the singularities 
of the type~$I_{i,j}$ for $i\leq k$ and~$j\leq l$. These strata are
also invariant with respect to the action of the group~$\Bbb C^*$. 
Taking the quotient modulo this action, we obtain the weighted
projective space
$\Bbb CP^{(k,l)}=\Bbb CP_{l,l,2l,\dots,(k-1)l,k,k,2k,\dots,(l-1)k}$,
and all the cohomology classes considered below belong to the cohomology
of this space.

Singularities of the type~$I_{i,j}$ 
have support on the subvariety~$\Delta$ of singular points of the 
fibers of the projection~$p$. 
In the coordinates chosen above this subvariety
is given by the equations $x=y=0$. Let us add to the classes
$\Psi,\Sigma\in H^2(\Bbb CP^{(k,l)})$ the class
$\Delta\in H^4(\Bbb CP^{(k,l)})$. We also consider the normal bundle
to~$\Delta$ and denote by~$N$ its first Chern class.
The class~$\Delta$ can be interpreted as the second Chern class
of this normal bundle, and the classes~$N$ and~$\Delta$ 
can be considered as (commuting) operators acting on~$\Delta$. 
Even less formally, one may treat~$N$ as a second cohomology 
class~$N\in H^2(\Bbb CP^{(k,l)})$ (which makes sense only 
when intersected with~$\Delta$). We are interested in expressions
for singularities of the type~$A_i$ in the unfoldings of singularities
of the type~$I_{k,l}$ in terms of the classes $\Sigma,\Psi,\Delta,N$.

\proclaim{Theorem 2.2}
There are expressions for the classes~$[A_i]$
in the unfoldings of singularities~$I_{k,l}$
as quasihomogeneous polynomials of degree~$i$ in
the variables $\Sigma,\Psi,\Delta,N$, such that{\rm:

(i)} these polynomials are universal, i.e., 
they depend on~$i$ but not on~$k$ and~$l${\rm;

(ii)} the polynomial corresponding to the class~$[A_i]$
can be represented as a sum of two polynomials,
the first summand coinciding with the polynomial~$P_i(\Sigma,\Psi)$
considered above, and the second summand having the support
on the stratum~$\Delta$ of double points{\rm;

(iii)} each polynomial supported on~$\Delta$
can be represented as a polynomial in~$\Delta,N$ and~$\Psi$.
\endproclaim

The first two statements follow from Thom's theorem
and Corollary~3.4 below, while the third one is a corollary
of the obvious relation $\Sigma\Delta=\Psi\Delta$.

Denote the polynomial corresponding to the class~$[A_i]$ by~$R_i$.
It has the form $R_i=P_i+Q_{i-2}\Delta$, where the polynomial~$P_i$
is already known and~$Q_{i-2}=Q_{i-2}(\Delta,N,\Psi)$
is a quasihomogeneous polynomial of degree~$i-2$.
We did not manage to find an explicit formula for these polynomials,
however, they can be computed using the indeterminate coefficients
method. Let us show how this method works for small~$i$.

Similarly to the above situation, each element in~$H^*(\Bbb CP^{(k,l)})$
can be expressed as a polynomial in~$t=\frac{\Psi}{kl}$\,. 
Knowing the weights of quasihomogenity, we conclude,
in particular, that
$$
\Psi=\Sigma=klt,\qquad N=(k+l)t,\qquad \Delta=klt^2.
$$

Let us find, using this representation, the coefficients
of the polynomials~$Q_0$ and~$Q_1$. We know that $Q_0=a$
for some constant~$a$ and that the stratum~$A_2$
in the unfolding of the singularity~$I_{1,1}$ is empty.
For $k=l=1$ we have $\Psi=\Sigma=t$, \ $\Delta=t^2$, whence
$$
R_2=\Sigma(2\Sigma-\Psi)+a\Delta=t^2+at^2=0.
$$
Therefore, $a=-1$, i.e.,
$$
R_2=\Sigma(2\Sigma-\Psi)-\Delta.
$$

Similarly, $Q_1=aN+b\Psi$ for some constants~$a$ and~$b$.
Using the fact that the stratum~$A_3$ is empty in the
unfoldings of each of the singularities~$I_{1,1}$
and~$I_{2,1}$, we obtain the following system of linear equations:
$$
\align
1+2a+b&=0,
\\
8+6a+4b&=0
\endalign
$$
(the homogeneity allowed us do divide both equations
by~$t^3$), whence $a=\nomathbreak 2$, \ $b=-5$, i.e.,
$$
Q_1=2N-5\Psi.
$$
Similar calculations yield
$$
\align
Q_2&=-(6N^2-15N\Psi+15\Psi^2-8\Delta),
\\
Q_3&=24N^3-62N^2\Psi+63N\Psi^2-35\Psi^3-60N\Delta+84\Psi\Delta,
\\
Q_4&=-(120N^4-322N^3\Psi+343N^2\Psi^2-196N\Psi^3+70\Psi^4
\\
&\qquad-
432N^2\Delta+812N\Psi\Delta-469\Psi^2\Delta+180\Delta^2).
\endalign
$$

Note that the number of indeterminate coefficients
in the expansion of~$A_i$ is $[i^2/4]$, 
where the square brackets denote the integral part
of a number, and this value exactly coincides with the
number of singularities of the type $I_{k,l}$
whose miniversal unfoldings do not contain the stratum~$A_i$ 
because of the dimensional reasons.

\subhead 2.3. Residual polynomials for multisingularities\endsubhead
Fix a partition $\alpha=1^{m_1}2^{m_2}\dots$
(where only finitely many exponents are nonzero).
Consider the set of points~$y$ in the space~$\Bbb C\times\Bbb C^m$
in the upper right corner of the commutative triangle~\thetag{2.1}
possessing the following properties:

-- the fiber~$p^{-1}(q(y))$ of the projection~$p$ is nonsingular;

-- the restriction of~$f$ to this fiber has, among the preimages of~$y$,
exactly~$m_1$ instances of the singularity~$A_1$, exactly~$m_2$
instances of the singularity~$A_2$ and so on (in addition to singularities
of the prescribed types, the point~$y$ may have also nonsingular
preimages, whose number is determined by the degree of~$f$).

Denote the class of the closure of this set times the order
$m_1!\,m_2!\,\dots$ of the automorphism group of the partition~$\alpha$
by $A_{\alpha}(f)\subset\Bbb C\times\Bbb C^m(f)=A_{1^{m_1}2^{m_2}\dots}(f)$.
(The multiplication by the order of the automorphism group
allows one to get rid of fractions in the further explicit
calculations of the classes.)

Were the pushforwards $f_*R_k$ of the subvarieties~$R_k=A_k$
intersecting transversally, the classes $A_{1^{m_1}2^{m_2}\dots}(f)$
would be described by the coefficients of
$\frac{t_1^{m_1}}{m_1!}\,\frac{t_2^{m_2}}{m_2!}\,\dots$
in the exponent of the generating function $f_*(R_1t_1+R_2t_2+\dots)$.
(The factorials $m_i!$ in the denominators reflect the possibility
to permute the preimages with the same singularity types.)
The nontransversality of the intersections makes
it necessary to introduce some correction terms, 
the residual polynomials~\cite{10}.

Consider the exponential generating function
$$
\Cal A (f)=\sum\bigl[A_{1^{m_1}2^{m_2}\dots}(f)\bigr]
\frac{t_1^{m_1}}{m_1!}\,\frac{t_2^{m_2}}{m_2!}\,\dots,
$$
where the summation is carried over all partitions 
$\alpha=1^{m_1}2^{m_2}\dots$, including the empty partition,
and the square brackets denote the cohomology classes
Poincar\'e dual to the corresponding subvarieties.
Let us introduce also the exponential generating function
in the infinite set of variables $t_1,t_2,\dots$,
with formal coefficients~$R_{\alpha}$, where~$\alpha$
takes values in the set of all nonzero partitions:
$$
\Cal R(t_1,t_2,\dots)=\sum R_{1^{m_1}2^{m_2}\dots}
\frac{t_1^{m_1}}{m_1!}\,\frac{t_2^{m_2}}{m_2!}\dots\,.
$$

\proclaim{Proposition 2.3}
There are weighted homogeneous polynomials~$R_{\alpha}$ in
the classes $\Sigma,\Psi,\Delta,N$ such that after substituting them
for the coefficients in the generating series~$\Cal R$
the cohomology classes of the strata~$A_{1^{m_1}2^{m_2}\dots}(f)$
coincide with the coefficients of the monomials
$\frac{t_1^{m_1}}{m_1!}\,\frac{t_2^{m_2}}{m_2!}\,\dots$
in the exponent $\exp(f_*\Cal R)${\rm:}
$$
\Cal A (f)=\exp\bigl(f_*\Cal R(t_1,t_2,\dots)\bigr).
$$
\endproclaim

Here the symbol~$f_*$ denotes the {\it Gysin homomorphism},
that is, the direct image homomorphism in the cohomology
of the quotient spaces modulo the actions of the group
$\Bbb C^*$. This mapping is a homomorphism of the additive, not
of the multiplicative structure. The polynomials~$R_{\alpha}$
introduced in this way are called the {\it residual polynomials
for multisingularities}. Of course, there is nothing strange
in the fact that the polynomials~$R_{k^1}$
coincide with the polynomials~$R_k$ introduced in Sec.\,2.2  
(indeed, the corresponding points have only one singular preimage).

\demo{Proof}
The most essential part of the Proposition states that
the correction term $f_*(R_{1^{m_1}2^{m_2}\dots})$
in the expression for the class $[A_{\alpha}]$ 
is linear with respect to the pushforwarded classes.
The main observation allowing one to prove this statement consists
in the fact that the universal expression for the multisingularities
classes remains valid even in the case of disconnected
underlying curves of meromorphic functions.

Let us apply this remark to a generic family of curves
$$
\alignedat2
&X&\overset{f}\to\longrightarrow\quad&\quad Y
\\
&\quad{}_{p}&\searrow\qquad\swarrow&{}_{q}\quad
\\
&&B\qquad &
\endalignedat
\tag 2.2
$$
Consider the mapping $f^{(d)}\:X^{(d)}\to Y$, where $X^{(d)}$
denotes the space fibered over~$B$, which is a disjoint union
of~$d$ copies of~$X$, \
$X^{(d)}=X_1\sqcup\dots\sqcup X_d$, \  $X_i=X$,
and~$f^{(d)}$ is the mapping whose restriction to the $i$~th
copy $X_i$ coincides with~$f$. To be more precise,
we take for the restriction of~$f^{(d)}$
to~$X_i$ a slightly perturbed copy of~$f$, which
is necessary for making~$f^{(d)}$ general.

Since the expression for classes of multisingularities
is universal, it can be applied to each of
the mappings~$f^{(d)}$ constructed above. It is easy to see that
$$
\Cal A (f^{(d)})=\bigl(\Cal A (f)\bigr)^d.
$$
Indeed, each partition~$\alpha$ corresponding to a multisingularity
of~$f^{(d)}$ splits into~$d$ subpartitions each
corresponding to a multisingularity of one of the components
of~$f^{(d)}$. The combinatorial coefficients of these partitions
are exactly the ones arising in the process of raising to the power~$d$.
We conclude that
$$
\log\Cal A (f^{(d)})=d\log\Cal A (f).
$$

On the other hand, it is obvious that the pushforward
of any monomial in the basic classes under the mapping~$f^{(d)}$ 
coincides with the $d$-fold pushforward of this monomial under~$f$.
Therefore, none of the coefficients 
in the generating function $\log\Cal A (f)$
contains a product of the pushforwards of basic classes
(the degree of the contribution of such coefficients would 
be greater than~$d$), whence all the coefficients are linear.
The Proposition is proved.
\enddemo

To compute the residual polynomials, we must know how
the homomorphism~$f_*$ acts on the cohomology classes 
we are interested in. The target of~$f$ also is 
a weighted homogeneous space, and its cohomology ring 
is generated by a hyperplane class, which we also denote by~$t$.
Hence $f_*$ has the form $f_*\:t^i\mapsto ct^i$
for some constant~$c$. This constant can be nothing but
the constant $f_*1$, that is, the degree of~$f$.
If~$f$ is a miniversal deformation of the singularity~$A_k$, 
then its degree equals~$k+1$.

Let us start with computing the polynomials~$R_{\alpha}$
for standard unfoldings of the singularities of the type~$A_k$
(i.e., in the presence of singularities of the type~$A$ solely).
These polynomials depend only on the classes~$\Sigma$
and~$\Psi$. They can be computed in the following way.

After substituting $\Sigma=kt$, \ $\Psi=(k+1)t$
into the generating function~$\Cal R$ we must obtain
$$
\split
&\exp\bigl((k+1)\Cal R(t_1,t_2,\dots)\bigr)=
1+(k+1)_2\frac{t_1}{1!}\,t+\biggl((k+1)_3\frac{t_2}{1!}+
(k+1)_4\frac{t_1^2}{2!}\biggr)t^2
\\
&\qquad\qquad+
\biggl((k+1)_4\frac{t_3}{1!}+(k+1)_5\frac{t_1}{1!}\,\frac{t_2}{1!}+
(k+1)_6\frac{t_1^3}{3!}\biggr)t^3+\dots\,.
\endsplit
\tag 2.3
$$

Indeed, the coefficient of the monomial
$\frac{t_1^{i_1}}{i_1!}\,\frac{t_2^{i_2}}{i_2!}\,\dots$ 
in the exponent must have the form
$t^{i_1+2i_2+3i_3+\dots}$ with a coefficient that is
a polynomial in~$k$ vanishing for all~$k$
up to the codimension of the corresponding stratum,
$k=-1,0,1,\dots,2i_1+3i_2+\dots-2$.

Note that the right-hand side of Eq.~\thetag{2.3}
can be conveniently rewritten in the form\footnote{We are
obliged to D.~Zagier for this remark.}
$$
\exp\bigl((k+1)\Cal R(t_1,t_2,\dots)\bigr)=
e^{\tau_1\tau+\tau_2\tau^2+\tau_3\tau^3+\dots}s^{k+1}\bigr|_{s=1}.
$$
Here $s$ is an additional variable commuting with all~$t_i$
and~$t$, and we use the notation
$$
\tau=t\frac{d}{ds}\,,\qquad \tau_i=t_i\frac{d}{ds}\,.
$$

Now it is easy to compute the residual polynomials~$R_{\alpha}$ 
by taking the logarithm
of the right-hand side of Eq.~\thetag{2.3},
dividing it by~$(k+1)$ and making the substitution
$t=\Psi-\Sigma$, \ $k=\Sigma/(\Psi-\Sigma)$.
For small codimensions we obtain
$$
\align
R_{1^2}&=-2\Sigma(3\Sigma-\Psi),
\\
R_{1^12^1}&=-6\Sigma(2\Sigma-\Psi)^2,
\\
R_{1^3}&=8\Sigma(15\Sigma^2-13\Sigma\Psi+3\Psi^2),
\\
R_{1^13^1}&=-4\Sigma(5\Sigma-3\Psi)(3\Sigma-2\Psi)(2\Sigma-\Psi),
\\
R_{2^2}&=-3\Sigma(2\Sigma-\Psi)(20\Sigma^2-25\Sigma\Psi+8\Psi^2),
\\
R_{1^22^1}&=24\Sigma(2\Sigma-\Psi)(15\Sigma^2-17\Sigma\Psi+5\Psi^2),
\\
R_{1^4}&=-48\Sigma(105\Sigma^3-160\Sigma^2\Psi+84\Sigma\Psi^2-15\Psi^3).
\endalign
$$

The $\Delta$-part of the residual polynomials in the
unfoldings of the singularities~$I_{k,l}$
can be computed recursively using the indeterminate coefficients
method, as above.

Let us compute, for example, the polynomial~$R_{1^2}$.
Due to the weighted homogeneity it has the form
$$
R_{1^2}=2\Sigma\Psi-6\Sigma^2+a\Delta
$$
for some constant~$a$ (we know already the
$\Delta$-free part), and the stratum ~$A_{1^2}$
can be represented as the unknown coefficient of $t_1^2/2!$
in $\exp(f_*\Cal R)$:
$$
[A_{1^2}]=(f_*R_{1^1})^2+f_*R_{1^2}=
(f_*\Sigma)^2+f_*(2\Sigma\Psi-6\Sigma^2+a\Delta).
$$

The degree of the unfolding of the singularity~$I_{k,l}$ is~$k+l$.
Substituting the expressions for the classes~$\Psi,\Delta,\Sigma$
in terms of the class~$t$ in the unfolding of the singularity~$I_{1,1}$
and making use of the fact that the class~$A_{1^2}$
is empty in this unfolding, we obtain the linear equation
$$
2^2+4-12+2a=0,
$$
whence $a=2$, i.e.,
$$
R_{1^2}=2(\Sigma\Psi-3\Sigma^2+\Delta).
$$

Similar calculations for residual polynomials in
higher codimensions lead to the following answers:
$$
\align
R_{1^12^1}&=-6\Sigma(2\Sigma-\Psi)^2+6\Delta(3\Psi-N_1),
\\
R_{1^3}&=8\Sigma(15\Sigma^2-13\Sigma\Psi+3\Psi^2)-8\Delta(10\Psi-3N),
\\
R_{1^13^1}&=-4\Sigma(5\Sigma-3\Psi)(3\Sigma-2\Psi)(2\Sigma-\Psi)
\\
&\qquad+
4\Delta(20\Psi^2-17N\Psi+6N^2-8\Delta),
\\
R_{2^2}&=-3\Sigma(2\Sigma-\Psi)(20\Sigma^2-25\Sigma\Psi+8\Psi^2)
\\
&\qquad+
3\Delta(25\Psi^2-21N\Psi+8N^2-12\Delta),
\\
R_{1^22^1}&=24\Sigma(2\Sigma-\Psi)(15\Sigma^2-17\Sigma\Psi+5\Psi^2)
\\
&\qquad-
24\Delta(20\Psi^2-15N\Psi+5N^2-7\Delta),
\\
R_{1^4}&=-48\Sigma(105\Sigma^3-160\Sigma^2\Psi+84\Sigma\Psi^2-15\Psi^3)
\\
&\qquad+
48\Delta(70\Psi^2-48N\Psi+15N^2-21\Delta).
\endalign
$$

\subhead 2.4. Residual polynomials for multimultisingularities\endsubhead
Now we are interested in subvarieties of the space in the lower corner of the
triangle~\thetag{2.1} consisting of points~$b$ such that
the restriction of~$f$ to the fiber~$p^{-1}(b)$ 
has several degenerate multisingularities, with prescribed 
reduced partitions~$\alpha_1^{m_1},\dts\alpha_c^{m_c}$ 
(as well as the necessary number of nondegenerate singularities).
The exponents in the notation indicate the numbers
of corresponding partitions. Denote by
$A_{\alpha_1^{m_1},\dts\alpha_c^{m_c}}(f)\subset\Bbb C^m$ 
the closure of the set of such points multiplied,
similarly to the multisingularity case, by the product
of the orders of the automorphism groups of the partitions~$\alpha_i$ 
and by the order of the automorphism group of the set
of the partitions~$\alpha_i$.

Our calculations show that the following statement, to which we 
did not find a formal proof yet, holds.

\proclaim{Conjecture 2.4}
Under the assumptions stated in the beginning of~{\S\,\rm2},
each class of multisingularities can be given as 
a universal polynomial combination of the images
of monomials in the basic classes under the Gysin homomorphism~$p_*$.
\endproclaim

Similarly to the multisingularities case, the very existence
of universal polynomials implies that they must have a rather special
form and be represented in terms of residual polynomials
in the basic classes. 

Consider the infinite set of variables 
$s_{1^2},s_{2^1},s_{1^2,1^2},s_{1^2,2^1},
s_{2^1,2^1},\dots$ indexed by sets of partitions~$\alpha$ 
(say, lexicographically ordered). Define a multiplication
on this set by setting the product of two variables 
equal to the sum of variables of the same kind,
whose indices are all sets of partitions that can be obtained
by uniting some pairwise distinct partitions in
the index of the second variable with some pairwise distinct partitions
in the index of the first variable. For example,
$$
s_{\alpha,\beta}s_{\gamma,\delta}=
s_{\alpha,\beta,\gamma,\delta}+s_{\alpha\gamma,\beta,\delta}+
s_{\alpha,\beta\gamma,\delta}+s_{\alpha\delta,\beta,\gamma}+
s_{\alpha,\beta\delta,\gamma}+s_{\alpha\gamma,\beta\delta}+
s_{\alpha\delta,\beta\gamma}.
$$
In particular, each polynomial in the variables~$s_{\alpha_1,\alpha_2,\dots}$
coincides with a linear polynomial.

Let us introduce the {\it generating function for multimultisingularities\/}
by the formula
$$
\frak A(s_{1^2},s_{2^1},\dots)=
\sum\frac{\bigl[A_{\alpha_1^{m_1},\dts\alpha_c^{m_c}}\bigr]}
{\prod \Aut(\alpha_i)}\,
\frac{s_{\alpha_1^{m_1},\dts\alpha_c^{m_c}}}{m_1!\ldots m_c!}\,.
$$
Consider also the exponential generating function 
with the formal coefficients
$T_{\alpha_1^{m_1},\dts\alpha_c^{m_c}}$,
indexed by the multipartitions:
$$
\frak R(s_{1^2},s_{2^1},\dots)=\sum
\frac{T_{\alpha_1^{m_1},\dts\alpha_c^{m_c}}}{\prod \Aut(\alpha_i)}\,
\frac{s_{\alpha_1^{m_1},\dts\alpha_c^{m_c}}}{m_1!\ldots m_c!}\,.
$$

\proclaim{Proposition 2.5}
Under Conjecture~{\rm2.4}, if we take for~$T_{\alpha}$ 
the polynomials representing the classes~$A_{\alpha}$, 
then there are weighted homogeneous polynomials
$R_{\alpha_1^{m_1},\dts\alpha_c^{m_c}}$, for $c\geq 2$,
in the classes $\Sigma,\Psi,\Delta,N$ 
such that after substituting their pushforwards
$f_*R_{\alpha_1,\dts\alpha_c}$ for the coefficients
$T_{\alpha_1^{m_1},\dts\alpha_c^{m_c}}$
in the above generating function~$\frak R$
the cohomology classes of the strata $A_{\alpha_1^{m_1},\dts\alpha_c^{m_c}}$
coincide with the coefficients of
$\dfrac{{s}_{\alpha_1^{m_1},\dts\alpha_c^{m_c}}}{m_1!\ldots m_c!}$ in
$\exp(q_*\frak R)$, where the exponent is taken with respect to
the multiplication of the variables~$s\phantom{\big|}$
described above.
\endproclaim

This statement is proved exactly in the same way 
as Proposition~2.3 by considering the mappings~$f^{(d)}$.

As above, the polynomials~$R_{\alpha_1,\dots,\alpha_c}$
can be computed inductively by means of the indeterminate
coefficients method. Calculations make use of the fact
that the restriction of~$p$ to~$\Sigma$ is a proper mapping,
and all the cohomology classes we are interested in
are supported on~$\Sigma$ and hence the homomorphism~$p_*$
expressed in terms of the basic class~$t$ is nothing but 
multiplication by the degree of the restriction
$p|_\Sigma$. The latter degree coincides with the
number of the critical points of~$f$ on a general fiber
and is equal to~$k+1$ for the singularity~$A_k$ 
and~$k+l$ for the singularity~$I_{k,l}$.
Here are the first nontrivial results of the calculations:
$$
\align
R_{2^1,2^1}&=-2\Sigma(2\Sigma-\Psi)(5\Sigma-3\Psi)+\Delta(14\Psi-5N),
\\
R_{2^1,1^2}&=6\Sigma(2\Sigma-\Psi)(5\Sigma-2\Psi)-6\Delta(7\Psi-2N),
\\
R_{1^2,1^2}&=-6\Sigma(30\Sigma^2-21\Psi\Sigma+4\Psi^2)+
2\Delta(57\Psi-14N).
\endalign
$$

\head
\S\,3. Relative Chern classes of the universal mapping
\endhead

In this section we show that the relative Chern classes
of the universal mapping indeed do belong to the subring
in~$H^*(X_{g,n})$ generated by the basic classes.
This means that the classes $[A_k]$ as well as the residual
polynomials for multisingularities and multimultisingularities
also belong to this subring. The main tool in the calculation
of relative Chern classes of~$f_{g,n}$ is the Grothendieck--Riemann--Roch
theorem. This theorem implies also nontrivial relations
for the Gysin pushforwards of polynomials in the basic classes.
We start with recalling the theorem in the specific situation
where we are going to use it. Our approach follows 
that of~\cite{3},~\cite{18}.

\subhead  3.1. Grothendieck--Riemann--Roch theorem
\endsubhead
The Grothendieck--Riemann--Roch (GRR) theorem 
relates the Chern character of a coherent sheaf over
a variety~$M$ with the Chern character of the direct image of this
sheaf under a mapping $P\:M\to A$. Suppose that both varieties~$M$
and~$A$ are compact and smooth.

In one of its forms (maybe, the most useful one)
the GRR formula looks like
$$
\ch(P_!\beta)=P_*\bigl(\ch(\beta)\td(P)^{-1}\bigr).
\tag 3.1
$$
Here $\beta$ is a coherent sheaf on~$M$,
\ $\ch$ denotes the Chern character, and~$\td$
denotes the Todd class of the morphism~$P$.

Let us give the definitions of the objects participating
in the GRR formula.

\subsubhead {\rm 3.1.1}. Chern character \endsubsubhead
The {\it Chern character\/} of the sheaf of sections
of a vector bundle, or simply the Chern character
of a vector bundle, is defined as the sum of the exponents
of its {\it Chern roots}, that is, the formal roots of 
its Chern polynomial,
$$
\ch(E)=\exp(e_1)+\cdots+\exp(e_r),
$$
where~$r$ denotes the rank of the vector bundle~$E$.
In particular, for a line bundle~$E$ we have
$\ch(E)=\exp\bigl(c_1(E)\bigr)$.

Being a symmetric function in the Chern roots, the Chern character
of a bundle can be expressed in terms of its characteristic classes:
$$
\split
\ch(E)&=r+c_1+\frac12\,(c_1^2-2c_2)+\frac16\,(c_1^3-3c_1c_2+3c_3)
\\
&\qquad+
\frac1{24}\,(c_1^4-4c_1^2c_2+4c_1c_3+2c_2^2-4c_4)
\\
&\qquad+
\frac1{120}\,(c_1^5-5c_1^3c_2+5c_1c_2^2+
5c_1^2c_3-5c_2c_3-5c_1c_4+5c_5)+\dots,
\endsplit
$$
where~$c_i$ denotes the~$i$~th symmetric function of the Chern roots.

Because of the properties of Chern roots, an exact sequence
of vector bundles 
$$
0\longrightarrow E'\longrightarrow E\longrightarrow E''\longrightarrow0
$$
leads to the {\it Whitney formula}
$$
\ch(E)=\ch(E')+\ch(E'').
$$

The Whitney formula allows one to extend the notion
of Chern character to arbitrary coherent sheaves.
A {\it coherent sheaf\/} on~$M$ is either a sheaf of sections 
of a vector bundle or, more generally, a sheaf
that can be inserted as the last nontrivial term
into an exact sequence of sheaves all whose terms but
the last one are sheaves of sections of vector bundles
(such an exact sequence is called the {\it resolvent\/}
for the given sheaf). The Chern character of a coherent sheaf
is defined as the alternated sum of the Chern characters
of the other terms in a resolvent.

For the structure sheaf of a subvariety~$\Sigma\subset M$ 
of codimension one in~$M$ given by an embedding
$i\:\Sigma\hookrightarrow M$,
the resolvent looks like follows:
$$
0\longrightarrow N_\Sigma^{\vee}\longrightarrow\Cal O_M
\longrightarrow i_*(\Cal O_\Sigma) \longrightarrow0,
$$
where $N_\Sigma$ is a line bundle over~$M$ such that~$\Sigma$ 
coincides with the set of zeros of its section.
In a tubular neighborhood of~$\Sigma$
the bundle~$N_\Sigma$ coincides with the natural
extension of the normal bundle to~$\Sigma$
in~$M$. The Whitney formula implies the identity
$$
\split
\ch(i_*(\Cal O_\Sigma))&=\ch(\Cal O_M)-\ch(N_\Sigma^{\vee})=
1-\exp(-\Sigma)
\\
&=
\frac{\bigl(1-\exp(-\Sigma)\bigr)}{\Sigma}\,\Sigma =
\biggl(\frac{\Sigma}{1-\exp(-\Sigma)}\biggr)^{-1}\Sigma,
\endsplit
$$
and the operator of intersection with the first Chern class
of the normal bundle can be replaced with~$\Sigma$,
since, by assumption, the variety 
$\Sigma$ coincides with the zero locus of a holomorphic section
of the normal bundle.

For a subvariety~$Z$ of codimension two the resolvent
proves to be more complicated, namely, it has the form
$$
0\longrightarrow \Lambda^2N_Z^{\vee}\longrightarrow N_Z^{\vee}
\longrightarrow\Cal O_M \longrightarrow i_*(\Cal O_Z) \longrightarrow 0,
\tag 3.2
$$
where~$\Lambda^2$ denotes the skew-symmetric square 
of a vector bundle. This sequence, called the 
{\it Koszul resolvent}, is written under the assumption
that~$Z$ is the zero locus of a globally defined section
of some rank two vector bundle $N_Z$ over~$M$.
We shall see that the final expression for
$\ch\bigl(i_*(\Cal O_Z)\bigr)$ 
is uniquely determined by the restriction of the bundle $N_Z$ to
$Z$, which coincides with the normal bundle to~$Z$. 
Moreover, one can show 
(see~\cite{3}) that the expression thus obtained remains
valid without assuming that the normal bundle to~$Z$
can be extended to the entire variety~$M$.
The Whitney formula implies
$$
\ch(i_*\Cal O_Z)=\ch(\Lambda^2N_Z^{\vee})+\ch(\Cal O_M)-\ch(N_Z^{\vee}).
$$

The rank of the virtual vector bundle $N_Z^{\vee}$
is~$2$, while its total Chern class is $c(N_Z^{\vee})=1-N_1+N_2$, where
$N_i=c_i(N_Z)$, \ $i=1,2$. Therefore, its Chern character
has the form
$$
\split
\ch(N_Z^{\vee})&=2-N_1+\frac12\,(N_1^2-2N_2)+\frac16\,(3N_1N_2-N_1^3)+
\frac1{24}\,(N_1^4-4N_1^2N_2+2N_2^2)
\\
&\qquad-
\frac1{120}\,(N_1^5-5N_1^3N_2+5N_1N_2^2)+\dots\,.
\endsplit
$$

The vector bundle $\Lambda^2N_Z^{\vee}$ is a line bundle and
$$
c(\Lambda^2N_Z^{\vee})=1-N_1,
$$
whence
$$
\ch(\Lambda^2N_Z^{\vee})=\exp(-N_1).
$$
Combining these results we obtain
$$
\align
\ch(i_*\Cal O_Z)&=N_2-\frac12\,N_1N_2+\ldots
\\
&
=\biggl(1-\frac12\,N_1+\frac1{12}\,(N_1^2-N_2)
-\frac1{24}\,(N_1^3-N_1N_2)+\dots\biggr)Z,
\endalign
$$
since both~$Z$ and~$N_2$ are identified with the
class of the zero section of the normal bundle $N_Z$.

\subsubhead {\rm 3.1.2.} Todd class \endsubsubhead
We start with the definition of the Todd class of a vector bundle,
and then define the Todd class of a mapping.
The {\it Todd class of a vector bundle\/}~$E$
over~$M$ is another symmetric function of its Chern roots,
$$
\split
\td(E)&=\frac{e_1}{1-\exp(-e_1)}\dots\frac{e_r}{1-\exp(-e_r)}
=1+\frac12\,c_1+\frac1{12}\,(c_1^2+c_2)
\\
&\qquad
+\frac1{24}\,c_1c_2+\frac1{720}\,(-c_1^4+4c_1^2c_2+3c_2^2+c_1c_3-c_4)+\dots\,.
\endsplit
$$
For the above exact sequence of vector bundles we have
$$
\td(E)=\td(E')\td(E'').
$$

For any vector bundle~$E$ of rank~$r$ we have
$$
\sum_{i=0}^r(-1)^i\ch(\Lambda^iE^{\vee})=c_r(E)\td(E)^{-1}.
$$

The {\it Todd class of a variety\/} $M$
is the Todd class of its tangent bundle~$TM$.
The {\it Todd class of a proper mapping\/}
$F\:M\to N$ is the Todd class of the difference $F^*TN-TM$, i.e.,
$$
\td(F)=\td(F^*TN-TM)=\frac{\td F^*TN}{\td TM}\,.
$$

\subsubhead {\rm 3.1.3}. Direct image of a coherent sheaf \endsubsubhead
The higher direct image~$P_!(E)$ of a coherent sheaf~$E$ on~$M$
with respect to a mapping $P\:M\to A$ is the most complicated
object among those defined in the present section. 
Methods of its computing depend on the difference between 
the dimensions of the varieties~$A$ and~$M$. In particular,
if the dimensions of the varieties coincide, then the image
$P_!(E)$ coincides with the pushforward~$P_*(E)$.
Recall that by definition the {\it pushforward\/} $P_*E$
of a sheaf~$E$ over~$M$ with respect to a mapping
$P\:M\to A$ is the sheaf whose group of sections over
a sufficiently small open set $U\subset A$
coincides with the restriction of~$E$ to
$P^{-1}(U)$. If $\dim M=\dim A+1$ (this is exactly the situation
where we are going to apply the GRR theorem), then
$$
P_!(E)=R^0\bigl(P_*(E)\bigr)-R^1\bigl(P_*(E)\bigr),
$$
where $R^0\bigl(P_*(E)\bigr)=P_*(E)$
is the pushforward of~$E$ and $R^1\bigl(P_*(E)\bigr)$
is the first higher Grothendieck direct image which can be
computed by means of the Serre duality:
$$
R^1\bigl(P_*(E)\bigr)=\bigl(R^0\bigl(P_*(E^{\vee}\otimes
\omega)\bigr)\bigr)^{\vee}=
\bigl(P_*(E^{\vee}\otimes\omega)\bigr)^{\vee};
$$
here $\omega$ is the relative dualizing sheaf.

\subhead
3.2. Relative Chern classes of the universal mapping
\endsubhead
Now let us return to the commutative triangle~\thetag{1.2} 
of spaces and mappings. From now on we shall use
simplified notation omitting the indices
$g,n$: we denote the space $X_{g,n}$
by~$X$, the mapping $f_{g,n}$ by~$f$, and so on. 
Introduce in~$X$ (recall that~$X$ is the universal curve over
the projectivized Hurwitz space) the following cohomology classes,
which we call {\it basic classes}. Denote by~$\omega$
the relative dualizing sheaf of the mapping $p\:X\to B$, i.e.,
the unique holomorphic line bundle over~$X$ whose fiber 
over a smooth point of a fiber of~$p$ coincides with the cotangent line
to this fiber. The first Chern class $c_1(\omega)$ 
is an element of the second cohomology group~$H^2(X)$. 
Let $\Pi$ be the divisor of poles of the universal mapping~$f$,
that is, the locus of points in~$X$ where the restriction of~$f$
to the fiber has a pole. Recall that since each
of the spaces $X$, $Y$ and~$B$ is a space of orbits under an action
of the group~$\Bbb C^*$, a cohomology class
$c_1\bigl(\Cal O(1)\bigr)\in H^2$ is distinguished on each of them;
we denote these classes by~$\Psi_X$, $\Psi_Y$ and~$\psi$, respectively. 
Introduce also the notation
$\Sigma=c_1(\omega)+\Psi_X+2\Pi$ for the first relative Chern class of~$f$.

In addition to subvarieties and classes of codimension one
we shall also make use of the class $\Delta\subset X$
represented by the set of singular points of the fibers of~$p$.
The rank of the normal bundle to~$\Delta$
in~$X$ is~$2$. Locally this normal bundle admits a natural
splitting into the direct sum of two line bundles.

Indeed, each point in~$\Delta$ is a point of transversal intersection
of the two branches of the fiber at this point, and that are the
tangent lines to these branches that produce the splitting of the
normal bundle. (In the case of rational curves, $g=0$, 
this splitting extends to a global splitting since a double point
cuts the rational curve into two connected components; we shall not
make use of this construction, however). Denote by~$N$
the first Chern class of the normal bundle to~$\Delta\subset X$. 
This Chern class can be treated as an operator that can be applied
to~$\Delta$. Its degrees $N^i$ are operators of the same form.
The second Chern class of the normal bundle can be identified
with the class~$\Delta$ itself.

\proclaim{Lemma 3.1}
The basic classes satisfy the following relations {\rm(}in the 
cohomology $H^*(X)${\rm):}
$$
\Pi^2=-\Psi_X \Pi,\qquad\Sigma \Pi=0,\qquad\Pi \Delta=0,\qquad
\Sigma \Delta =\Psi_X \Delta.
$$
\endproclaim

\demo{Proof}
The identity $\Pi \Delta=0$ is a corollary of the fact that
the poles of a function in the Hurwitz space are smooth points of the
underlying curve, whence the subvarieties
$\Delta$ and~$\Pi$ are disjoint. Similarly, the identity
$\Sigma \Pi=0$ expresses the fact that all the poles of a function
are simple, and there is no ramification over infinity.

The dualizing sheaf $\omega$ is almost trivial over $\Delta$
(its pullback to the two-fold covering of $\Delta$ is trivial),
and it is isomorphic to the conormal line bundle
to the hypersurface $\Pi$ at points of this hypersurface.
These properties imply the identities
$c_1(\omega) \Delta=0$ and~$\Pi^2=-c_1(\omega) \Pi$.
Substituting $c_1(\omega)=\Sigma-
\Psi_X-2\Pi$ into these equations we obtain, respectively,
$$
\gather
0=(\Sigma-\Psi_X-2\Pi) \Delta=\Sigma \Delta-\Psi_X \Delta,
\\
\Pi^2=-(\Sigma-\Psi_X-2\Pi) \Pi=\Psi_X \Pi+2\Pi^2,
\endgather
$$
and the identities
$\Sigma \Delta =\Psi_X \Delta$ and~$\Pi^2=-\Psi_X \Pi$ follow.
Note that the latter can also be obtained by pulling back to~$X$
the identity $\Pi_Y^2=-\Psi_Y \Pi_Y$, which takes place on~$Y$
(see Sec.~3.3 below).
\enddemo

\proclaim{Corollary 3.2}
Each polynomial in the basic classes can be represented in the form
$$
P_1(\Psi_X) \Pi+P_2(\Psi_X,\Sigma)+ P_3(\Psi_X,N,\Delta) \Delta,
$$
where $P_1$, $P_2$, and $P_3$ are some polynomials.
\endproclaim

The main goal of the present section consists
in the proof of the following statement, from which
we shall deduce some relations for the pushforwards
of the basic classes.

\proclaim{Theorem 3.3}
The total relative Chern class of~$f$ is
$$
c(f)=(1+\Psi_X)\biggl(\frac1{1-\Sigma+\Psi_X}-
\frac\Delta{1+N+\Delta}\biggr).
\tag 3.3
$$
\endproclaim

The theorem immediately implies

\proclaim{Corollary 3.4}
If the mapping~$f$ can have singularities only of types~$A_k$
or $I_{k,l}$, when restricted to the fibers of~$p$,
then the residual polynomials~$R_{\alpha_1,\dots,\alpha_c}$ 
for multisingularities belong to the subring in~$H^*(X)$
generated by the basic classes $\Sigma,\Psi_X,\Delta,N$.
\endproclaim

\goodbreak
Since~$p$ is the composition of the mappings~$f$ and~$q$, \ $p=f\circ q$,
its total Chern class is the ratio of the pullback of the 
total relative Chern class of~$q$ and the total relative Chern class of~$f$:
$$
c(f)=\frac{f^*c(q)}{c(p)}\,.
$$
Therefore, in order to compute this class it suffices
to compute the classes~$c(p)$, $c(q)$,
as well as the pullback of the latter class to~$X$.
This is precisely what we are going to do now.

\subsubhead {\rm 3.2.1}. The class~$c(q)$ 
and its pullback~$f^*c(q)$ \endsubsubhead
The total space~$Y$ of the bundle~$q$ is the projectivization
of the vector bundle $E=\Bbb C\oplus\Cal O(1)$ of rank~$2$ over~$B$.
Denote by $\Cal T $ the tautological line bundle over~$Y$ and set
$\Pi_Y=c_1(\Cal T ^{\vee})=-c_1(\Cal T )$. This class
satisfies the relation $c_2\bigl(q^*(E)/\Cal T\bigr)=0$, i.e.,
$$
\Pi_Y^2+\Pi_Y \Psi_Y=0,
\tag 3.4
$$
where $\Psi_Y=q^*(\psi)=c_1(\Cal O(1))$
is the characteristic class of the action of the group~$\Bbb C^*$. 
The resulting relation describes the cohomology ring of~$Y$ 
as an algebra over $H^*(B)$,
and the action of the ring $H^*(B)$ 
is determined by the mapping~$q^*$.

The lines $\Bbb C\oplus\{0\}$ and~$\{0\}\oplus\Cal O(1)$
in the fibers of~$E$ correspond to the points
``0'' and ``$\infty$'' in the fibers of~$q$. These points
form hypersurfaces with the Poincar\'e dual classes
$$
c_1\bigl(\Hom\bigl(\Cal T,\Cal O(1)\bigr)\bigr)=\Pi_Y+\Psi_Y,\qquad
c_1\bigl(\Hom\bigl(\Cal T,\Bbb C)\bigr)=\Pi_Y
$$
respectively, hence the class $\Pi_Y$ can be interpreted
as the ``dual class to the set of infinities in~$Y$'',
and Eq.~\thetag{3.4} expresses the fact that the set of zeros
and the set of infinities are disjoint.

The total Chern class of~$q$ is inverse to the total
Chern class of the relative tangent bundle $\ker(q_*)$.
This bundle is given by the isomorphism
$\ker(q_*)\simeq\Hom\bigl(\Cal T,q^*(E)/\Cal T\bigr)$. Therefore,
$$
c(q)=\frac1{c\bigl(\Hom\bigl(\Cal T,q^*(E)/\Cal T\bigr)\bigr)}=
\frac1{1+\Psi_Y+2\Pi_Y}\,.
$$
Applying $f^*$ we obtain $f^*(\Psi_Y)=\Psi_X$ and~$f^*(\Pi_Y)=\Pi$
and, finally,
$$
f^*\bigl(c(q)\bigr)=\frac1{1+\Psi_X+2\Pi}\,.
$$

\subsubhead {\rm 3.2.2}. The total Chern class of~$p$ \endsubsubhead
The computation of the total Chern class of the mapping~$p$
is slightly more complicated. It can be computed using the
exact sequence
$$
0\longrightarrow
p^*T^{\vee} B\overset{p^{\vee}}\to\longrightarrow T^{\vee} X
\longrightarrow\omega\longrightarrow\Cal O_{\Delta}\longrightarrow0,
\tag 3.5
$$
where, recall,~$\omega$ denotes the relative dualizing sheaf
of~$p$. The subvariety $\Delta\subset X$ is exactly
the set of those points where the corank of $p^*$ 
fails to be one and is~$2$. The structure sheaf~$\Cal O_{\Delta}$
of the subvariety~$\Delta$ measures the difference between the relative
dualizing sheaf~$\omega$ and the sheaf of relative differentials.

Applying to the above exact sequence the Whitney formula, we obtain
$$
c(p^{\vee})=1-c_1(p)+c_2(p)-\ldots=
c(p^*T^{\vee} B-T^{\vee} X)=\frac{c(\Cal O_\Delta)}{c(\omega)}=
\frac{c(\Cal O_\Delta)}{1+c_1(\omega)}\,.
$$

To compute the Chern class~$c(\Cal O_\Delta)$
let us use the Koszul projective resolvent~\thetag{3.2}.
Applying the Whitney formula once again we have
 $$
c(\Cal O_{\Delta})=\frac{1-N}{1-N+\Delta}=1-\frac\Delta{1-N+\Delta}\,.
$$
This formula can be deduced as well
by applying the GRR formula to the embedding
$\Delta\hookrightarrow X$. Substituting it in the expression for
$c(p^{\vee})$, we obtain
$$
c(p)=\frac1{1-c_1(\omega)}\biggl(1-\frac\Delta{1+N+\Delta}\biggr).
$$

Now we can compute the total Chern class of~$f$.
Applying the relations in Lemma~3.1, we have
$$
\split
c(f)&=\frac{c(p)}{f^*\bigl(c(q)\bigr)}=\frac{1+\Psi_X+2\Pi}{1-c_1(\omega)}
\biggl(1-\frac\Delta{1+N+\Delta}\biggr)
\\
&=
\frac{1+\Psi_X+2\Pi}{1-\Sigma+\Psi_X+2\Pi}-
\frac{1+\Psi_X}{1+N+\Delta}\Delta
\\
&=
\biggl(\frac{1+\Psi_X}{1-\Sigma+\Psi_X}-
\frac{1+\Psi_X}{1+N+\Delta}\Delta\biggr)+
\biggl(\frac{1+\Psi_X+2\Pi}{1-\Sigma+\Psi_X+2\Pi}-
\frac{1+\Psi_X}{1-\Sigma+\Psi_X}\biggr).
\endsplit
$$
The first summand coincides with the right-hand side of Eq.~\thetag{3.3};
after reducing the second summand to the common denominator
its numerator becomes
$$
(1+\Psi_X+2\Pi)(1-\Sigma+\Psi_X)-(1+\Psi_X)(1-\Sigma+\Psi_X+2\Pi)=
-2\Pi \Sigma=0.
$$
Theorem~3.3 is proved.

\subhead 3.3. Pushforward homomorphisms\endsubhead
In order to be able to apply the universal formula
to classes of multimultisingularities, we must have 
a description of the pushforward homomorphisms $f_*$, $p_*$, $q_*$.
The following theorem states that the homomorphism~$p_*$ 
uniquely determines the two others.

\proclaim{Theorem 3.5}
For any set of cohomology classes $h_1,\dots,h_s$ on~$X$
satisfying the property $h_i \Pi=0$, \ $i=1,\dots,s$,
we have
$$
q_*\bigl(f_*(h_1)\dots f_*(h_s)\bigr)=\psi^{s-1}p_*(h_1)\dots p_*(h_s).
$$
\endproclaim

The assumption $h_i\Pi=0$ means geometrically that the class $h_i$
admits a representation by a cycle nonintersecting the pole divisor~$\Pi$. 
It is satisfied for all the residual polynomials of multisingularities.

\demo{Proof}
Since $q$ is the projectivization of a vector bundle,
the direct image $f_*(h)$ of an arbitrary class $h\in H^*(X)$
can be expressed in the form
$$
f_*(h)=q^*(a)+q^*(b) \Pi_Y
$$
for some classes $a$ and~$b$ on~$B$. The class $b$
can be computed by applying $q_*$ to both sides of this identity:
the identity $q_*\Pi_Y=1$ together with the projection formula yields
$$
b=q_*\bigl(f_*(h)\bigr)=p_*(h).
$$

The class~$a$ can be computed in a similar way:
multiplying both sides of the above equation by~$\Pi_Y+\Psi_Y$,
we obtain
$$
f_*\bigl((\Psi_X+\Pi) h\bigr)=q^*(a) (\Psi_Y+\Pi_Y)+
q^*(b) (\Psi_Y+\Pi_Y) \Pi_Y=q^*(\psi a)+q^*(a) \Pi_Y,
$$
and applying $q_*$ we get $a=p_*\bigl((\Psi_X+\Pi) h\bigr)$.
If the class~$h$ satisfies the assumption of the theorem,
$h \Pi=0$, then the answer can be simplified:
$$
f_*(h)=q^*\bigl(p_*(h)\bigr)(\Psi_Y+\Pi_Y).
$$
Therefore, under the assumptions of the theorem we have
$$
f_*(h_1)\dots f_*(h_s)=(\Psi_Y+\Pi_Y)^s
\prod_{i=1}^s q^*\bigl(p_*(h_i)\bigr)=
(\Psi_Y^s+\Psi_Y^{s-1}\Pi_Y)\prod_{i=1}^s q^*\bigl(p_*(h_i)\bigr).
$$
Applying $q_*$ and the projection formula, we obtain the desired
equation.
\enddemo

Hence to complete the calculation of the cohomology classes
of the strata in the Hurwitz space 
$P\overline{\Cal H}_{g,n}$ 
we must know how to compute the pushforward
$p_*(R)$ of any residual polynomial expressed in the form
$$
R=P_1(\Psi_X,\Sigma)+P_2(\Psi_X,N,\Delta) \Delta.
$$

Let us set
$$
\xi_k=p_*(\Sigma^{k+1}),\qquad \delta_{k,\ell}=p_*(N^k\Delta^{\ell+1}).
$$
There are no universal relations between the classes
$\xi_k$ and~$\delta_{k,\ell}$ that are valid for all values of~$g$ and~$n$. 
Of course, for any given~$g$ and~$n$ the dimension of the cohomology
space $H^*(P\overline{\Cal H}_{g,n})$ is finite, hence there
are such relations. Some of them are discussed in Sec.3.4 below.

\subhead 3.4. Relations on cohomology classes in~$B$\endsubhead
Let us apply the GRR formula to the projection
$p\:X\to B$ and to the relative dualizing sheaf~$\omega$.
We proceed similarly to~\cite{18}, where the same thing has been done
for the projection of the universal curve to the moduli space
of curves. We have
$$
\ch(p_!\omega)=p_*\bigl(\ch(\omega)\td(p)^{-1}\bigr).
$$
The sheaf $p_!\omega=R^0p_*\omega-R^1p_*\omega$ on the left-hand side
can be expressed in terms of the Hodge bundle~$\Lambda=\Lambda_{g,n}$ 
over~$B=P\overline{\Cal H}_{g,n}$ (recall that the fiber
of the Hodge bundle over a point $(C;f)$ is the space
of holomorphic $1$-forms over~$C$; hence the rank of the bundle is~$g$):
\ $R^0p_*\omega=p_*\omega=\Lambda$ and, by the Serre duality,
$R^1p_*\omega=p_*\Bbb C=\nomathbreak \Bbb C$. 
The classes on the right-hand side of the formula were computed above.
Using the relation $c_1(\omega)\Delta=0$, we obtain
$$
 \split
\ch(\omega)\td^{-1}(p)&=\ch(\omega)\td(\omega^{\vee}-\Cal O_\Delta^{\vee})=
e^{c_1(\omega)}\frac{-c_1(\omega)}{1-e^{c_1(\omega)}}\,
\td(-\Cal O_\Delta^{\vee})
\\
&=\frac{c_1(\omega)}{1-e^{-c_1(\omega)}}+\td(-\Cal O_\Delta^{\vee}),
\endsplit
$$
where the Todd class of the (virtual) sheaf $-\Cal O_\Delta^{\vee}$
is totally determined by its Chern class
$$
c(-\Cal O_\Delta^{\vee})=\frac{1+N+\Delta}{1+N}=1+\frac\Delta{1+N}\,.
$$

Substituting this result into the GRR formula, we obtain
$$
\split
&(g-1)+\ch_1(\Lambda)+\ch_2(\Lambda)+\ch_3(\Lambda)+
\ch_4(\Lambda)+\ch_5(\Lambda)+\dots
\\
&\qquad=
\frac12\,p_*\bigl(c_1(\omega)\bigr)+
\frac1{12}\,p_*\bigl(c^2_1(\omega)+\Delta\bigr)-
\frac1{720}\,p_*\bigl(c^4_1(\omega)+(N^2-3\Delta)\Delta\bigr)
\\
&\qquad\qquad+
\frac1{30240}\,p_*\bigl(c^6_1(\omega)+
(N^4-5N^2\Delta+5\Delta^2)\Delta\bigr)+\dots\,.
\endsplit
\tag 3.6
$$

This relation shows that the homogeneous components 
of the Chern character (whence, the Chern classes)
of the Hodge bundle can be expressed in terms of the pushforwards
of the basic classes.

Note that there are no terms of positive even degrees
on the right-hand side of Eq.~\thetag{3.6}.
This means that the even part of the character
$\ch(\Lambda)$  is trivial, i.e., $c(\Lambda+\Lambda^{\vee})=1$.
If $g=0$, then the Hodge bundle itself is
trivial, whence the left-hand side of Eq.~\thetag{3.6} 
vanishes at positive degrees. This property implies
nontrivial identities for the pushforwards of the basic classes:
$$
 p_*\bigl(c^2_1(\omega)+\Delta\bigr)=0,\qquad
 p_*\bigl(c^4_1(\omega)+(N^2-3\Delta)\Delta\bigr)=0,\qquad \dots\,.
$$

\head
\S\,4. Applying universal polynomials to studying
the stratification of Hurwitz spaces
\endhead

The knowledge of the residual polynomials for
singularities, multisingularities, and multimultisingularities
allows one to compute the cohomology classes of the strata 
in the Hurwitz spaces. Complete formulas of this kind
must take into account the contribution of nonisolated singularities
as well, which we did not manage to compute yet.
That is why we present the universal expressions only
for strata in the case of rational functions,
and only for small codimensions.

\subhead  4.1. Universal expressions for strata\endsubhead
Below we present the expressions for strata of codimensions~$1$
and~$2$ modulo classes supported on strata of nonisolated
singularities in terms of direct images of the basic classes
obtained on the base of the residual classes computed above
(we use the notation
$\xi_k=p_*(\Sigma^{k+1})$, \ $\xi_0=2n-2+2g$, \
$\delta_{k,l}=p_*(N^k\Delta^{l+1})$):
$$
\align
\sigma_{2^1}&=-\psi\xi_0+2\xi_1-\delta_{0,0},
\\
\sigma_{1^2}&=\frac12\,\psi\xi_0 (\xi_0+2)-3 \xi_1+\delta_{0,0},
\\
\sigma_{3^1}&=2 \xi_0 \psi^2 -7 \xi_1 \psi-5 \delta_{0,0} \psi+
6 \xi_2+2 \delta_{1,0},
\\
\sigma_{1^12^1}&=-\xi_0 (\xi_0+6) \psi^2+2 (\xi_0+12)\xi_1\psi+
(18-\xi_0) \delta_{0,0} \psi-24 \xi_2-6 \delta_{1,0},
\\
\sigma_{1^3}&=\frac16 \,\xi_0 (\xi_0^2+6 \xi_0+24) \psi^2
-\frac13 \,(9 \xi_0+52) \xi_1 \psi
\\
&\qquad+
\frac13\,(3\xi_0-40) \delta_{0,0} \psi+20 \xi_2+4 \delta_{1,0},
\\\vspace{3mm}
\sigma_{2^1,2^1}&=\frac12\, (\xi_0-6) \xi_0 \psi^2-(2 \xi_0  - 11)\xi_1\psi+
\frac12 \,\delta_{0,0}^2+2 \xi_1^2-10 \xi_2
\\
&\qquad+
(\xi_0+7) \delta_{0,0}\psi-2 \xi_1 \delta_{0,0}-\frac52\, \delta_{1,0},
\\
\sigma_{2^1,1^2}&=-\frac12\, \xi_0 (\xi_0^2-2 \xi_0-12) \psi^2+
(\xi_0^2+ \xi_0-27) \xi_1 \psi-\delta_{0,0}^2-6 \xi_1^2+30 \xi_2
\\
&\qquad-
\frac12\, (\xi_0^2+42)\delta_{0,0} \psi+5 \xi_1 \delta_{0,0}+6 \delta_{1,0},
\\
\sigma_{1^2,1^2}&=\frac18\, (\xi_0-4) \xi_0 (\xi_0^2+4 \xi_0+6) \psi^2-
\frac34\, (2 \xi_0^2-4 \xi_0-21) \xi_1 \psi+\frac12\, \delta_{0,0}^2
\\
&\qquad+
\frac92\, \xi_1^2-\frac{45}2\, \xi_2+
\frac14 \,(2 \xi_0^2-4 \xi_0+57) \delta_{0,0}\psi-
3 \xi_1 \delta_{0,0}-\frac72\, \delta_{1,0}.
\endalign
$$

Now let us give the expressions for the same strata 
in the case~$g=0$ simplified with the help of the identities
$$
\xi_0=2n-2,\qquad\xi_1=4 (n - 1) \psi-\delta_{0,0}
$$
deduced in Sec.~3.4 (note that in the case
$g=0$ the codimension of the locus of nonisolated singularities in~$B$
is~$3$, and hence they do not contribute to the above equations):
$$
\align
\sigma_{2^1}&=6(n-1) \psi-3 \delta_{0,0},
\\
\sigma_{1^2}&=2(n-6) (n-1) \psi+4 \delta_{0,0},
\allowdisplaybreak
\\
\sigma_{3^1}&= -24 (n-1) \psi^2+2 \delta_{0,0} \psi+6 \xi_2+2\delta_{1,0},
\\
\sigma_{1^12^1} &= 12 (n-1) (n+6) \psi^2-6 n \delta_{0,0}\psi-
24 \xi_2-6 \delta_{1,0},
\\
\sigma_{1^3} &= \frac43\, (n-1) (n^2-17 n-30) \psi^2+
4 (2 n-1) \delta_{0,0} \psi+20 \xi_2+4 \delta_{1,0},
\\\vspace{3mm}
\sigma_{2^1,2^1} &= 2 (n-1) (9 n+10) \psi^2-2 (9 n  -7) \delta_{0,0} \psi+
\frac92\, \delta_{0,0}^2  -10 \xi_2-\frac52\,\delta_{1,0},
\\
\sigma_{2^1,1^2} &= 12 (n-9) (n-1) n \psi^2-
6 (n^2-13 n+11) \delta_{0,0} \psi-12 \delta_{0,0}^2+30\xi_2+6\delta_{1,0},
\\
\sigma_{1^2,1^2} &= (n-1) (2 n^3-30 n^2+145 n-60) \psi^2+
\frac12\, (16 n^2-144 n+125) \delta_{0,0} \psi
\\
&\qquad+
8\delta_{0,0}^2-\frac{45}2\, \xi_2-\frac72\, \delta_{1,0}.
\endalign
$$

\subhead 4.2. Degrees of strata in genus zero\endsubhead
In the case $g=0$ the compactified Hurwitz space 
$P\overline{\Cal H}_{0,n}$ is fibered over the compactified
moduli space $\overline{\Cal M}_{0,n}$ 
of rational curves with marked points. The fiber of this
bundle is the projective space $PE$, where 
$E=L_1^{\vee}\oplus \cdots\oplus L_n^{\vee}$
is the Whitney sum of the tangent lines to the curve
at its marked points. This means that the cohomology
of the space $P\overline{\Cal H}_{0,n}$
treated as an algebra over $H^*(\overline{\Cal M}_{0,n})$
are generated by the class $\psi=c_1(\Cal O(1))$
subject to the relation
$$
\psi^n+c_1(E) \psi^{n-1}+\cdots+c_n(E)=0
$$
(cf.~\cite{16}). In particular, each class $\alpha$ in
$H^{2d}(P\overline{\Cal H}_{0,n})$ can be represented in the form
$$
\alpha=\pi^*(\eta_{d})+\pi^*(\eta_{d-1})\psi+
\pi^*(\eta_{d-2})\psi^2+\dots,
$$
where $\pi\:P\overline{\Cal H}_{0,n}\to\overline{\Cal M}_{0,n}$
is the natural projection, and the classes $\eta_i$ 
belong to the cohomology of the space
$\overline{\Cal M}_{0,n}$. The well-known relation
$\pi_*\psi^s=c_{s-n+1}(-E)$ (see, e.g.,~\cite{3}) 
implies that the degree of any class represented in this form
can be computed using the relation
$$
 \deg\alpha=\int_{P\overline{\Cal H}_{0,n}}
 \frac{\pi^*(\eta_{d})+\pi^*(\eta_{d-1})\psi+\dots}{1-\psi}=
 \int_{\overline{\Cal M}_{0,n}}
 \frac{\eta_{d}+\eta_{d-1}+\eta_{d-2}+\dots}{c(E)}\,.
$$
This relation reduces the computation of the degree
to a computation of some integrals over the moduli spaces
of rational curves with marked points, whose values
are well known (see, e.g., \cite{2} or~\cite{14}).
Below we represent the results of computations of 
the degrees of the pushforwarded basic classes up to codimension~$2$:
$$
\align
\deg(1) &= n^{n-3},
\\
\deg(\delta_{0,0}) &= \frac12\, (n-1) (n+6) n^{n-4},
\\
\deg(\xi_2) &= \frac13\, (n-1) (17 n^2-28 n+12) n^{n-5},
\\
\deg(\delta_{1,0}) &= -\frac16\, (n-1) (n^2+10 n-120) n^{n-5},
\\
\deg(\delta_{0,0}^2) &= \frac1{12}\, (n-1) (3 n^3+31 n^2+82 n-120) n^{n-5}.
\endalign
$$

Substituting these equations in the expressions for the strata
we obtain the following formulas for their degrees:
$$
\allowdisplaybreaks
\align
\deg(\sigma_{2^1})&=\frac92\, (n-2) (n-1) n^{n-4},
\\
\deg(\sigma_{1^2})&=2 (n-3) (n-2) (n-1) n^{n-4},
\\
\deg(\sigma_{3^1})&=\frac{32}3\, (n-3) (n-2) (n-1) n^{n-5},
\\
\deg(\sigma_{1^12^1}) &= 9 (n-4) (n-3) (n-2) (n-1) n^{n-5},
\\
\deg(\sigma_{1^3})&= \frac43\, (n-5) (n-4) (n-3)(n-2) (n-1) n^{n-5},
\\\vspace{3mm}
\deg(\sigma_{2^1,2^1})&=\frac38\, (n-2)(n-1) (27 n^2-137 n+180) n^{n-5},
\\
\deg(\sigma_{2^1,1^2})&=3 (n-3) (n-2)(n-1)(3n^2-15 n+20) n^{n-5},\\
\deg(\sigma_{1^2,1^2}) &= (n-3) (n-2) (n-1)(2n^3-16 n^2+43 n-40) n^{n-5}.
\endalign
$$

Now Theorem~1.1 produces the following formulas for the Hurwitz numbers:
$$
\allowdisplaybreaks
\gather
{
\align
 h_{2^1} &= \frac92 \frac{(2 {n-4})!}{(n-3)!}\, n^{n-5},
 \\
 h_{1^2} &= 2 \frac{(2 n-4)!}{(n-4)!}\, n^{n-5},
 \\
 h_{3^1} &= \frac{32}3 \frac{(2 n-5)!}{(n-4)!}\, n^{n-6},
 \\
 h_{1^12^1} &= 9\, \frac{(2 n-5)!}{(n-5)!}\, n^{n-6},
 \\
 h_{1^3} &= \frac43 \frac{(2 n-5)!}{(n-6)!}\, n^{n-6},
\endalign}
\\\vspace{3mm}
{
\align
 h_{2^1,2^1}&=\frac34\, (27 n^2-137 n+180)\frac{(2 n-6)!}{(n-3)!}\, n^{n-6},
 \\
 h_{2^1,1^2} &= 3 (3 n^2-15 n+20)\frac{(2 n-6)!}{(n-4)!}\, n^{n-6},
 \\
 h_{1^2,1^2}&=2 (2 n^3-16 n^2+43 n-40)\frac{(2 n-6)!}{(n-4)!}\, n^{n-6};
\endalign}
\endgather
$$
the first group of formulas here refers to the case
of a single degenerate ramification, and these formulas
are special cases of the Hurwitz formula. 

Another approach to the calculation of the degrees of the strata
consists in the direct application of the Hurwitz formula
(or, for an arbitrary genus, the formula from~\cite{2}
which is, however, less explicit) for Hurwitz numbers of multisingularities.
Theorem~1.1 allows one to find the degree of each 
stratum of the form~$\sigma_{\alpha}$. 
In the rational case $g=0$ knowing the decompositions of the strata
as polynomials in the pushforwards of basic classes
we can compute the degrees of at least some of these direct images.
In our case it is easy to compute the degrees
$\deg(1)$, \ $\deg(\delta_{0,0})$, \
$\deg(\delta_{1,0})$, \ $\deg(\xi_2)$. 
Note, however, that knowing these degrees is insufficient for
computing the degrees of the strata of multimultisingularities:
their decomposition includes the class~$\delta_{0,0}^2$.
The degree of the latter can be computed, for example,
by using the degree of the stratum~$\sigma_{2^1,2^1}$, 
which is known from~\cite{22}.

\subhead 4.3. On nonisolated singularities\endsubhead
Although we do not yet have a complete understanding of what
is going on when nonisolated singularities are added to the picture,
a few words about them are in order. Consider the first nontrivial
contribution of nonisolated singularities in the case
of rational functions, i.e., the classes of codimension~$3$.
Denote by $\sigma^{\exp}_{\alpha_1,\alpha_2,\dots}$ 
the expected decompositions for the corresponding classes
obtained by the universal formulas under the assumption that
there are no singularities of types other than~$A_k$
and~$I_{k,l}$. They differ from the correct decompositions of the classes
$\sigma_{\alpha_1,\alpha_2,\dots}$ by terms supported on the locus
$I_{\infty}\subset B$ of functions with nonisolated
singularities. Since the correction terms have the same dimension
as the subvariety $I_{\infty}$, they must be proportional to the
fundamental class of this subvariety, with a constant proportionality
coefficient. It is natural to assume that this coefficient is universal,
that is, it depends neither on~$n$, nor on the choice
of an irreducible component of the subvariety $I_{\infty}$. 
The following ``unexpected'' observation serves a numerical
confirmation of this assumption: the degree of the difference
$$
\deg\bigl(\sigma_{\alpha_1,\alpha_2,\dots}-
\sigma^{\exp}_{\alpha_1,\alpha_2,\dots}\bigr)
$$
is proportional to the degree of the stratum $I_{\infty}$,
$$
\deg(I_{\infty})=\frac18\,(n-1)(n^3+11n^2+34n-120)n^{n-5},
$$
with integer coefficients. The degrees of the classes
$\sigma^{\exp}_{\alpha_1,\alpha_2,\dots}$ can be computed 
by means of the methods described in the present paper, while
the correct degree of the stratum $\sigma_{\alpha_1,\alpha_2,\dots}$ 
can be deduced from the Hurwitz formula.
Let us present the conjectural contribution of nonisolated
singularities to some strata of codimension three:
$$
\align
\sigma_{4^1}&=\sigma^{\exp}_{4^1}+5 I_{\infty},
\\
\sigma_{1^13^1}&=\sigma^{\exp}_{1^13^1}-16 I_{\infty},
\\
\sigma_{2^2}&=\sigma^{\exp}_{2^2}-9 I_{\infty},
\\
\sigma_{1^22^1}&=\sigma^{\exp}_{1^22^1}+36 I_{\infty},
\\
\sigma_{1^4}&=\sigma^{\exp}_{1^4}-16 I_{\infty}.
\endalign
$$
\nobreak

\enddocument